\documentclass[11pt]{amsart}
\usepackage{amssymb,latexsym, amsmath, amsxtra}
\usepackage[dvips]{graphics}
\usepackage{epsfig}

\theoremstyle{plain}
\newtheorem{theorem}{Theorem}[section]

\newtheorem{proposition}[theorem]{Proposition}
\theoremstyle{definition}

\theoremstyle{remark}

\numberwithin{equation}{section}

\begin{document}


\title{Triple correlation of the Riemann zeros}

\abstract We use the conjecture of Conrey, Farmer and Zirnbauer
for averages of ratios of the Riemann zeta function \cite{kn:cfz2}
to calculate all the lower order terms of the triple correlation
function of the Riemann zeros.  A previous approach was suggested
by Bogomolny and Keating \cite{kn:bk96} taking inspiration from
semi-classical methods.  At that point they did not write out the
answer explicitly, so we do that here, illustrating that by our
method all the lower order terms down to the constant can be
calculated rigourously if one assumes the ratios conjecture of
Conrey, Farmer and Zirnbauer.  Bogomolny and Keating
\cite{kn:bk06} returned to their previous results simultaneously
with this current work, and have written out the full expression.
The result presented in this paper agrees precisely with their
formula, as well as with our numerical computations, which we
include here.

We also include an alternate proof of the triple correlation of
eigenvalues from random $U(N)$ matrices which follows a nearly
identical method to that for the Riemann zeros, but is based on
the theorem for averages of ratios of characteristic polynomials
\cite{kn:cfz1,kn:cfs05}.
\endabstract

\author {J.B. Conrey}
\address{American Institute of Mathematics,
360 Portage Ave, Palo Alto, CA 94306} \address{School of
Mathematics, University of Bristol, Bristol, BS8 1TW, United
Kingdom} \email{conrey@aimath.org}

\author{N.C. Snaith}
\address{School of Mathematics,
University of Bristol, Bristol, BS8 1TW, United Kingdom}
\email{N.C.Snaith@bris.ac.uk}

\thanks{
Research of the first author supported by the American Institute
of Mathematics. The second author was supported by an EPSRC
Advanced Research Fellowship. Both authors have been supported by
a Focused Research Group grant (0244660) from the National Science
Foundation} \maketitle

\tableofcontents

\section{Introduction}
\label{sec:intro}

In 1973 Montgomery \cite{kn:mont73} proved the following result,
assuming the Riemann Hypothesis, on the two-point correlation of
the zeros of the Riemann zeta function:
\begin{eqnarray}
\label{eq:mont} &&\sum_{\gamma_1,\gamma_2\in [0,T]}
w(\gamma_1-\gamma_2)
f\left(\frac{\log T}{2\pi}(\gamma_1-\gamma_2)\right)\\
&& \qquad\qquad = \frac{T\log T}{2\pi} \left( f(0) + \int
_{-\infty}^{\infty} f(u) \left[ 1-\left(\frac{\sin (\pi u)}{\pi u}
\right)^2\right] du +o(1)\right)\nonumber
\end{eqnarray}
for suitably decaying functions $f$ with Fourier transform
supported in $[-1,1]$ and weight $w(x)=\frac{4}{4+u^2}$. He
conjectured that (\ref{eq:mont}) would in fact hold for any test
function $f$.

In 1994 Hejal \cite{kn:hejhal94} proved a similar result for the
triple correlation of Riemann zeros:
\begin{eqnarray}
&&\sum_{{\gamma_1,\gamma_2,\gamma_3\in[T,2T]}\atop {\rm distinct}}
w(\gamma_1,\gamma_2,\gamma_3) f\bigg(\frac{\log T}{2\pi}(\gamma_1
-\gamma_2),\frac{\log T}{2\pi}(\gamma_1-\gamma_3)\bigg) \\
&&\qquad=\frac{T\log T}{2\pi} \bigg( \int_{-\infty}^{\infty} \int
_{-\infty}^{\infty}f(u,v)\left| \begin{array}{ccc} 1& S(u) & S(v)
\\ S(u)&1&S(u-v)\\S(v)&S(u-v)&1\end{array}\right| du\;dv \nonumber\\
&&\qquad\qquad\qquad+o(1)\bigg)\nonumber
\end{eqnarray}
with weight $w(x_1,x_2,x_3)= \prod _{j<k} \exp \big[-\frac{1}{6}
(x_j-x_k)^2\big]$ and the Fourier transform of the continuous,
suitably decaying test function $f$ is supported on the hexagon
with successive vertices $(1,0)$, $(0,1)$, $(-1,1)$, $(-1,0)$,
$(0,-1)$ and $(1,-1)$.  Here $S(x)=\frac{\sin(\pi x)} {\pi x}$.

This was extended to the $n$-point correlation function and to
more general $L$-functions by Rudnick and Sarnak \cite{kn:rudsar}
in 1996.

\begin{figure}[htbp]
  \begin{center}
    \includegraphics[scale=1.0]
    {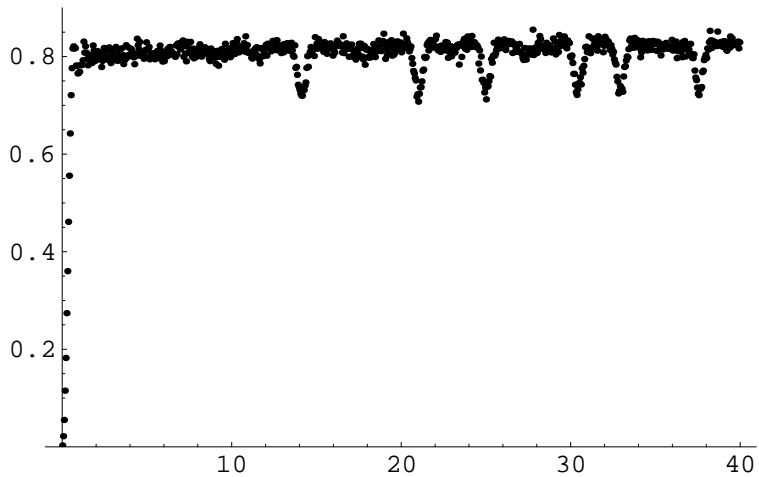}
        \caption{The two-point correlation function of the Riemann zeta function
calculated from the first 100 000 zeros.  The correlation distance
is plotted along the x-axis.  Note the significant dip at the
first few Riemann zeros: 14.13, 21.02, 25.01, $\ldots$.}
    \label{fig:2point}
  \end{center}
\end{figure}

These results encompass the rigorous work on the subject, but are
limited by two things.  Firstly, the support of the Fourier
transform of the test function is always confined to a restricted
range.  Secondly, only the asymptotic for large $T$ is found. This
second point is understandable, as the goal was to show that this
limiting form was the same as that for the $n$-point correlation
function of eigenvalues from large-dimensional matrices from the
GUE ensemble of random matrix theory (see \cite{kn:conrey1} or
\cite{kn:keasna03} for review articles on the connection between
random matrix theory and number theory).  This aim was duly
achieved, but there is clearly interest in the lower-order terms,
as Bogomolny and Keating's early results \cite{kn:bk96} showed
that in the two-point correlation function of the Riemann zeros
one sees sensitivity to the positions of the low Riemann zeros
themselves - something that clearly does not happen in random
matrix theory at any order. After it was predicted by Bogomolny
and Keating, a numerical illustration of this for the two-point
correlation function was first shown in \cite{kn:berrykeating99},
where Berry and Keating also fully explain  a similar phenomenon
in the number variance statistic first observed by Berry in 1988
\cite{kn:berry88}.  A numerical plot of the two-point correlation
function calculated using the first 100 000 zeros of the Riemann
zeta function is shown in Figure \ref{fig:2point}.
The x-axis is the correlation distance between pairs of zeros in
unscaled units, showing the distinctive dip at each Riemann zero.
The plot is a histogram of the number of pairs of zeros with a
given separation distance, and the y-axis is divided by
$T\Big(\frac{\log\tfrac{T}{2\pi}}{2\pi}\Big)^2$, where $T\sim 75
000$ is the height of the 100 000th zero.

Meanwhile, alongside this rigorous work in the number theory
community, physicists using semiclassical techniques applied in
the field of quantum chaos, treated the Riemann zeta function as a
model system (the prime numbers playing the role of periodic
orbits) and so studied the correlation of Riemann zeros in
analogous ways to those in which they would study correlation of
energy levels in their more standard physical systems.  The first
step in this direction was the derivation by Keating
\cite{kn:keating93} of the limiting form of the two-point
correlation function of the Riemann zeros.  This is a heuristic
calculation and it relies on a conjecture by Hardy and Littlewood
\cite{kn:harlitt18} for the behaviour of correlations between
prime numbers, but it has the advantage that there are no
restrictions on the support of a test function.  Using analogous
methods this result was then extended in two papers by Bogomolny
and Keating \cite{kn:bogkea95,kn:bogkea96} to obtain the limiting
form of the $n$-point correlation function.

The first result on lower-order terms of the correlations of the
Riemann zeros was also by Bogomolny and Keating \cite{kn:bk96},
but using a different heuristic inspired more directly from
semiclassical methods.  Here they truncate the Euler product for
the Riemann zeta function at primes less than $\log
\tfrac{T}{2\pi}$ (in semiclassical language this means taking
periodic orbits up to the Heisenberg time) and define a new set of
zeros from the resulting approximation to the staircase function
of the Riemann zeros (the function that increases by one at the
position of each zero).  It is the two-point correlation of this
new set of zeros that miraculously gives all the significant
lower-order terms of the analogous statistic for the Riemann
zeros.  Using the same method, an expression is also given
\cite{kn:bk96} in semiclassical language for all lower-order terms
of the three-point correlation function, and it is this which
could be turned into a detailed formula including all terms
calculated in the present paper, but the authors did not publish
it explicitly at that time.

Bogomolny and Keating also did further work, largely unpublished,
obtaining all the lower-order terms for the 2-, 3- and 4-point
correlations both by extending their Hardy-Littlewood method (see
\cite{kn:keating99} for some details of the two-point correlation
function) and by the method mentioned in the previous paragraph,
as well as two other heuristic methods (see \cite{kn:bog00}).

Recently \cite{kn:bk06} they have written out the lower-order
terms for the three point correlation function of the Riemann
zeros in full detail and these agree with the results presented in
this paper.

\section{The Riemann zeros}
\label{sect:rzf}

\subsection{Results}

Using the ratios conjecture of Conrey, Farmer and Zirnbauer
\cite{kn:cfz1,kn:cfz2} we obtain the following

\begin{theorem}  \label{theo:triple} Assuming the ratios conjecture and summing over distinct zeros of the Riemann zeta
function:
\begin{eqnarray}
\label{eq:triple} &&\sum_{0<\gamma_1\neq\gamma_2\neq\gamma_3<T}
f(\gamma_1-\gamma_2, \gamma_1-\gamma_3)
=\frac{1}{(2\pi)^3}\int_{-T}^T\int_{-T}^T
f(v_1,v_2)\nonumber \\
&& \qquad \qquad\times\Bigg[\int_0^T \log ^3 \frac{u}{2\pi} du +
I(iv_1,iv_2;0) + I(0,iv_1;-iv_2) + I(0,iv_2;-iv_1)\nonumber
\\
&&\qquad\qquad\qquad\qquad +
I(-iv_1,-iv_2;0)+I(0,-iv_2;iv_1)+I(0,-iv_1;iv_2)\nonumber
\\&&\qquad\qquad\qquad+I_1(0;iv_2)+I_1(0;iv_1)+I_1(-iv_2;iv_1) +I_1(-iv_2;0)
+I_1(-iv_1;iv_2)\nonumber
\\
&&\qquad\qquad\qquad\qquad+I_1(-iv_1;0)\Bigg]dv_1dv_2+O(T^{\epsilon}),
\end{eqnarray}
where the integrals in $v_1$ and $v_2$ are to be interpreted as
principal value integrals and

\begin{eqnarray}\label{eq:Iaab}
I(\alpha_1,\alpha_2;\beta)&:=& \int_0^T\frac{\zeta'}{\zeta}
(\tfrac{1}{2} +it+\alpha_1)\; \frac{\zeta'}{\zeta} (\tfrac{1}{2}
+it+\alpha_2)\;\frac{\zeta'}{\zeta} (\tfrac{1}{2} -it+\beta)dt
\nonumber \\
&=&\int_0^T
Q(\alpha_1+\beta,\alpha_2+\beta)+\Big(\frac{t}{2\pi}\Big)^{-\alpha_1-\beta}
\zeta(1-\alpha_1-\beta)
\zeta(1+\alpha_1+\beta)\nonumber \\
&&\qquad\times\bigg(A(\alpha_1+\beta)\Big(\frac{\zeta'}{\zeta}(1+\alpha_2-\alpha_1)-\frac
{\zeta'}{\zeta}(1+\alpha_2+\beta)\Big)+P(\alpha_1+\beta,\alpha_2+\beta)\bigg)\nonumber
\\
&& +\Big(\frac{t}{2\pi}\Big)^{-\alpha_2-\beta}
\zeta(1-\alpha_2-\beta)
\zeta(1+\alpha_2+\beta)\nonumber \\
&&\qquad\times\bigg(A(\alpha_2+\beta)\Big(\frac{\zeta'}{\zeta}(1+\alpha_1-\alpha_2)-\frac
{\zeta'}{\zeta}(1+\alpha_1+\beta)\Big)+P(\alpha_2+\beta,\alpha_1+\beta)\bigg)dt\\
&&\qquad\qquad\qquad\qquad+O(T^{1/2+\epsilon}),\nonumber
\end{eqnarray}
\begin{eqnarray}
\label{eq:I1ab} I_1(\alpha;\beta)&:=& \int_0^T \log
\tfrac{t}{2\pi}\; \frac{\zeta'}{\zeta} (\tfrac{1}{2}
+it+\alpha)\;\frac{\zeta'}{\zeta} (\tfrac{1}{2} -it+\beta)dt
\nonumber \\
&=&\int_0^T\log \tfrac{t}{2\pi}\bigg(\Big(\frac{\zeta'}{\zeta}\Big)'(1+\alpha+\beta) \nonumber\\
&&\qquad+\Big(\frac{t}{2\pi} \Big)^{-\alpha-\beta}
\zeta(1+\alpha+\beta)\zeta(1-\alpha-\beta) A(\alpha+\beta) -
B(\alpha+\beta) \bigg)dt \\
&&\qquad\qquad\qquad\qquad+O(T^{1/2+\epsilon})\nonumber
\end{eqnarray} and
\begin{eqnarray}
\label{eq:Aprime} A(x)&=&\prod_p
\frac{(1-\tfrac{1}{p^{1+x}})(1-\tfrac{2}{p}+\tfrac{1}{p^{1+x}})}
{(1-\tfrac{1}{p})^2},
\end{eqnarray}
\begin{eqnarray}
\label{eq:Bprime} B(x)&=&\sum_p\left(\frac{\log
p}{p^{1+x}-1}\right)^2,
\end{eqnarray}
\begin{eqnarray}
\label{eq:Qprime} Q(x,y)&=&-\sum_p \frac{\log^3 p}
{p^{2+x+y}(1-\tfrac{1}{p^{1+x}})( 1-\tfrac{1}{p^{1+y}})}
\end{eqnarray}
and
\begin{eqnarray}
\label{eq:Pprime}
P(x,y)&=&\prod_p\frac{\Big(1-\frac{1}{p^{1+x}}\Big)\Big(1-\frac{2}{p}+\frac{1}{p^{1+x}}\Big)
} {\Big(1-\frac{1}{p}\Big)^2} \\
&&\qquad\times \sum_p\frac{\Big(1-\frac{1}{p^x}\Big)
\Big(1-\frac{1}{p^x}-\frac{1}{p^y} +\frac{1}{p^{1+y}}\Big) \log
\tfrac{1}{p}} {\Big(\frac{1}{1-p^{1-x+y}}\Big)
\Big(1-\frac{1}{p^{1+y}}\Big) \Big(
1-\frac{2}{p}+\frac{1}{p^{1+x}}\Big)p^{2-x+y}}.\nonumber
\end{eqnarray}
\end{theorem}

\begin{figure}[htbp]
  \begin{center}
    \includegraphics[scale=1.0]
    {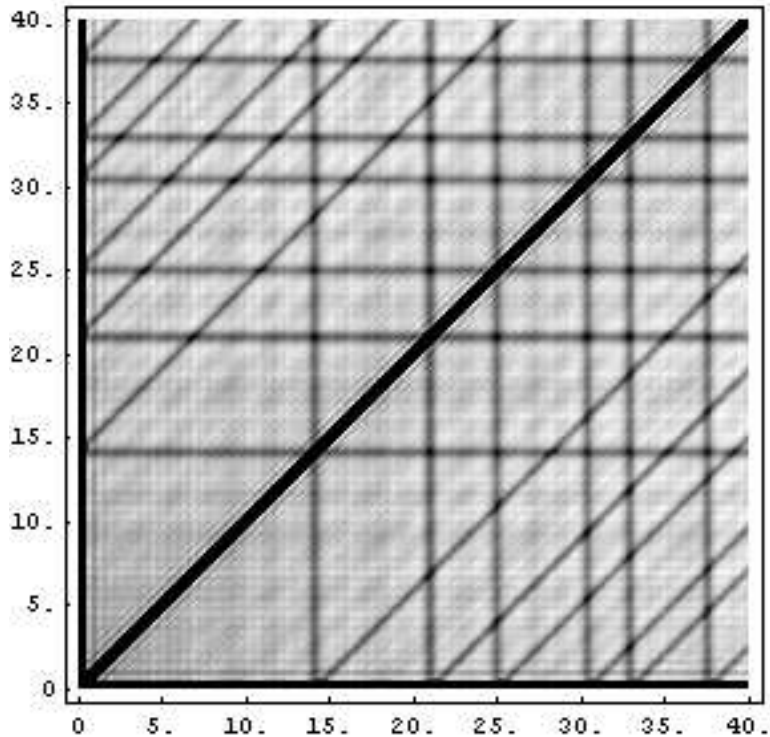}
    \caption{The triple correlation of the Riemann zeros: we plot
    the quantity in square brackets from (\ref{eq:triple}).  The x
    and y axes are $v_1$ and $v_2$ and the density plot is lighter were
    the function has a higher value and is darker for smaller
    values. Note
the horizontal, vertical and diagonal lines occurring at the
Riemann zeros: 14.13, 21.02, 25.01, $\ldots$}
    \label{fig:triptheory}
  \end{center}
\end{figure}

\begin{figure}[htbp]
  \begin{center}
    \includegraphics[scale=1.0]
    {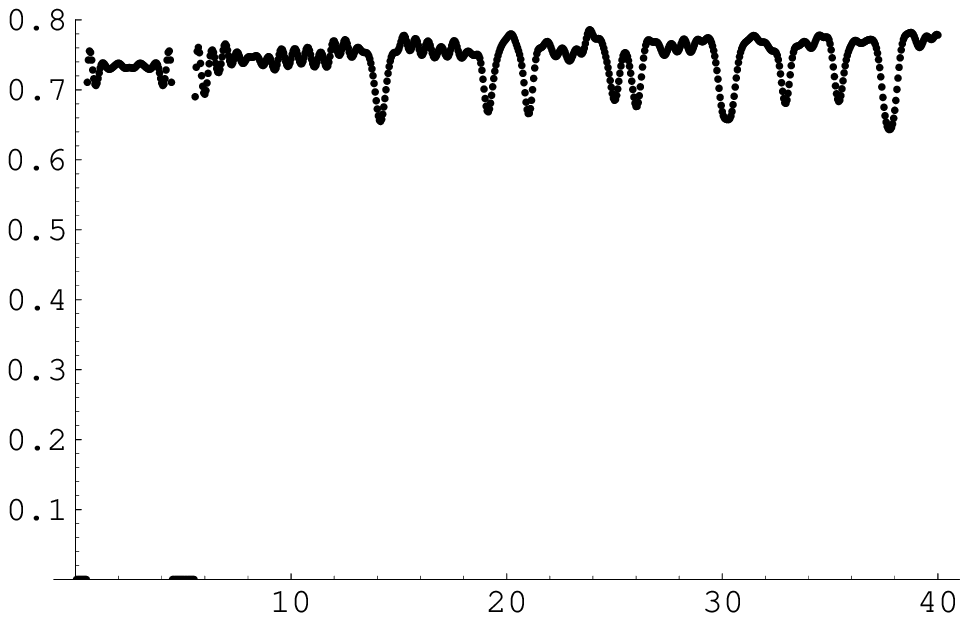}
    \caption{A horizontal cross-section of Figure
\ref{fig:triptheory} at height 5 on the y-axis to illustrate the
height of that plot. }
    \label{fig:crossect}
  \end{center}
\end{figure}

The expression (\ref{eq:triple}) is plotted in Figure
\ref{fig:triptheory}, with
$f(v_1,v_2)=\delta(v_1-x)\delta(v_2-y)$. The x- and y- axes are
unscaled, but (\ref{eq:triple}) is divided by
$T\Big(\frac{\log\tfrac{T}{2\pi}}{2\pi}\Big)^3$.  The density plot
is light for large values and dark where (\ref{eq:triple}) is
small.  Note the horizontal, vertical and diagonal lines occurring
at the Riemann zeros, caused by terms like
$\frac{\zeta'}{\zeta}(1+ix)$, $\frac{\zeta'}{\zeta}(1+iy-ix)$,
etc.   Numerical computation of (\ref{eq:triple}) breaks down near
the x- and y- axes and on the diagonal because of the principal
value integration, so the plot has been set to zero in these
regions. The plot could be completed with a more careful expansion
of the formula around $x=0$, $y=0$ and $x=y$, but this would not
be particularly edifying.  The result would be extremely similar
to the random matrix limit shown in Figure \ref{fig:triprmt}.  The
maximum height on the contour plot in Figure \ref{fig:triptheory}
is about 0.799.  An idea of the height of the plot can be seen in
Figure \ref{fig:crossect}, which is a horizontal cross-section of
Figure \ref{fig:triptheory} at height 5 on the y-axis.

\begin{figure}[htbp]
  \begin{center}
    \includegraphics[scale=1.0]
    {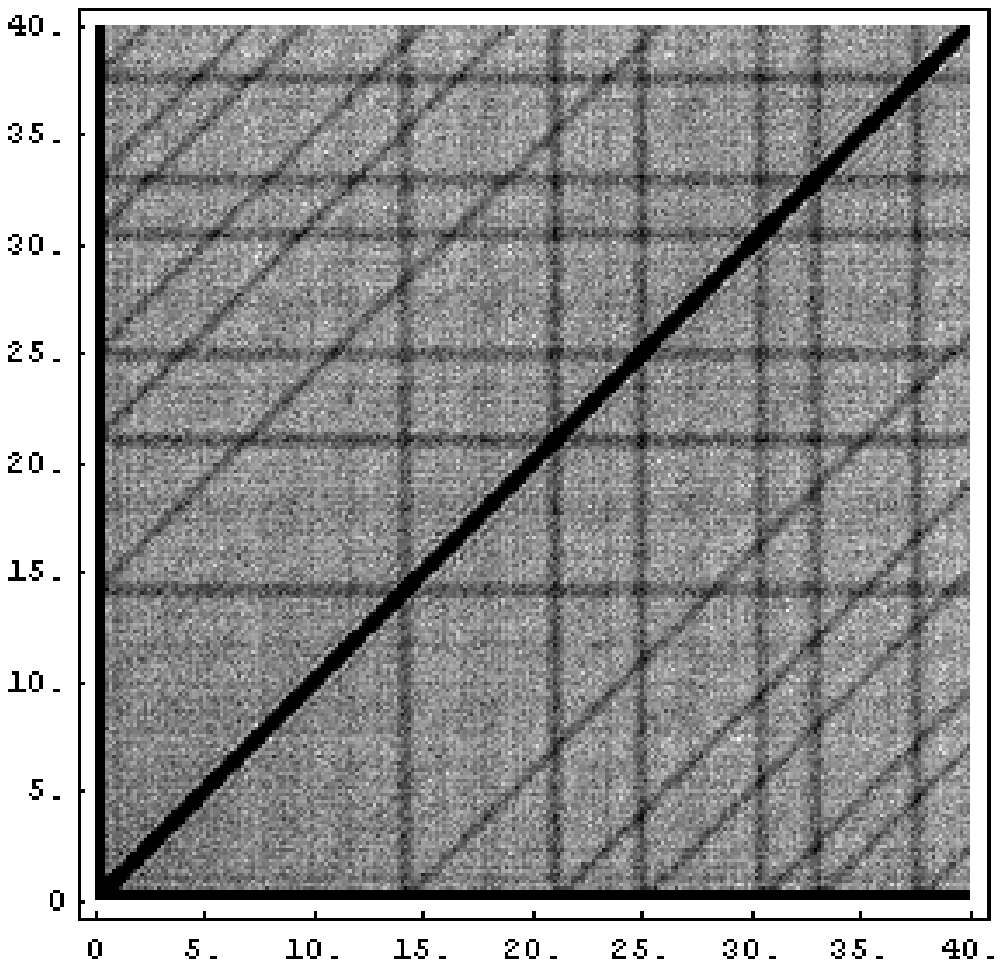}
    \caption{The numerically calculated triple correlation of
the Riemann zeros.  The distribution of triplets of zeros is
depicted by plotting  the distance between the first and second
zero in the triplet on the horizontal axis against the distance
between the first and third zero on the vertical. Higher
occurrences appear lighter grey and where no triplets fall the
plot is black.}
    \label{fig:tripnum}
  \end{center}
\end{figure}

\begin{figure}[htbp]
  \begin{center}
    \includegraphics[scale=1.0]
    {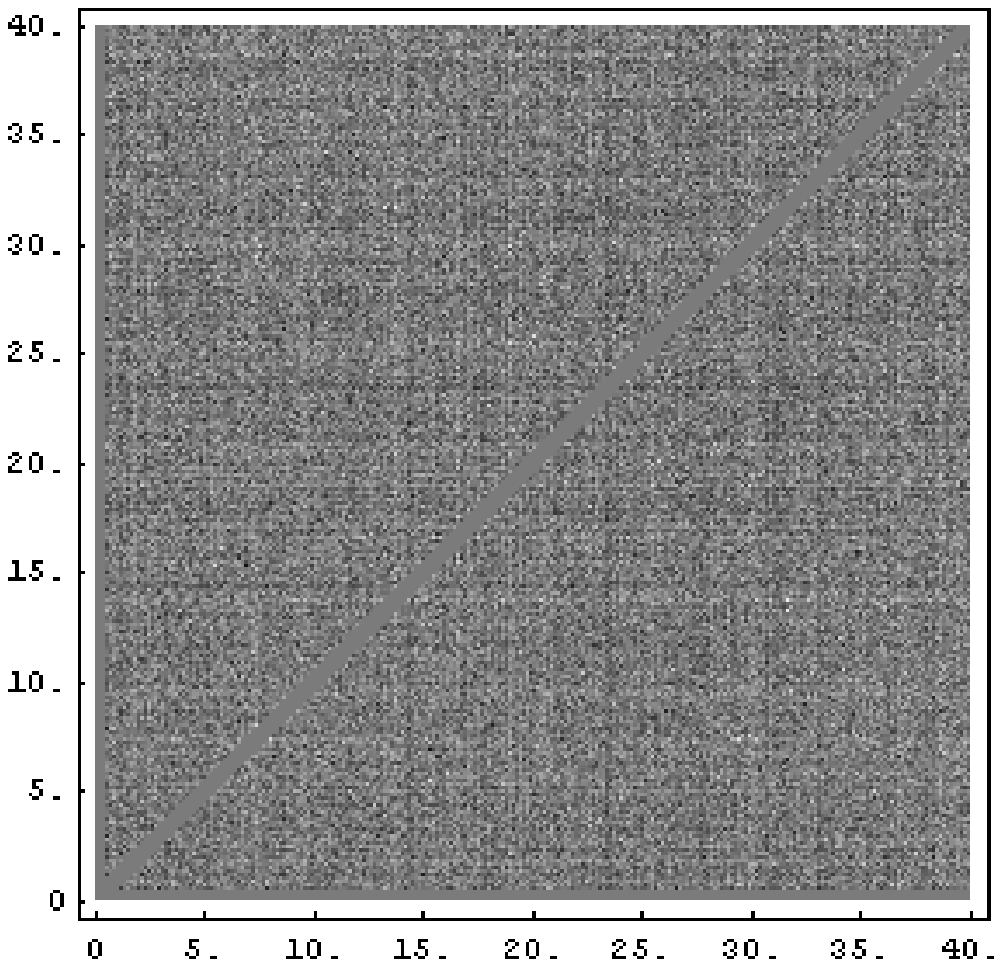}
    \caption{The difference between Figures
\ref{fig:tripnum} and \ref{fig:triptheory}.}
    \label{fig:tripdiff}
  \end{center}
\end{figure}

In Figure \ref{fig:tripnum} the numerical triple correlation,
using the first 100 000 zeros, is plotted, again scaled as above.
The difference between this and Figure \ref{fig:triptheory} is
shown in Figure \ref{fig:tripdiff}, where the maximum height of
the plot is about 0.164.  Compare this with the maximum height of
Figure \ref{fig:tripnum}, which is about 0.923.  The maximum value
of the difference plot may seem rather large, as we expect an
error of $T^{-1/2+\epsilon}$, but this is probably due to the
relatively small value of $T$ used for these plots. For $T\sim
75000$, $T^{-1/2}$ is around 0.003, but for $T$ values of this
size powers of $\log T$ can make a big difference.  The mean value
of the points on the difference plot (Figure \ref{fig:tripdiff})
is -0.00127, and the standard deviation is about 0.03, which gives
a better idea of the spread of the points. We also note that the
mean of the absolute value of Figure \ref{fig:tripdiff} is about
0.0257.

For ease of comparison with Theorem \ref{theo:triple}, we state
here the similar random matrix result which is derived in detail
in Section \ref{sect:rmtratios}.  The identical structure of
Theorem \ref{theo:triple} and Theorem \ref{theo:tripleRMT} is
apparent if we recall the equivalence $N=\log \frac{t}{2\pi}$ (see
for example \cite{kn:keasna00}) and reduce (\ref{eq:tripleRMT})
from three to two variables by performing a simple translation
such as $\theta_2\rightarrow \theta_2+\theta_1$ and
$\theta_3\rightarrow \theta_3+\theta_1$ and noting from
(\ref{eq:IaabRMT}) and (\ref{eq:I1abRMT}) that
$J(\alpha,\beta;\gamma)=J(\alpha+t,\beta+t;\gamma-t)$ and
$J(\alpha;\beta)=J(\alpha+t;\beta-t)$.  In comparing
$J(\alpha,\beta;\gamma)$ with $I(\alpha,\beta;\gamma)$ and
$J(\alpha;\beta)$ with $I_1(\alpha;\beta)$, $z(x)=(1-e^{-x})^{-1}$
plays the part of $\zeta(1+x)$, as is always the case in moment
and ratios conjectures.   Using the ratios theorem of Conrey,
Farmer and Zirnbauer \cite{kn:cfz1,kn:cfz2} we obtain the
following

\begin{theorem}  \label{theo:tripleRMT} With the star indicating a sum over distinct
eigenvalues we have:
\begin{eqnarray}
\label{eq:tripleRMT}&&
T_3(f):=\int_{U(N)}\sideset{}{^*}\sum_{j_1,j_2,j_3}f(\theta_{j_1},\theta_{j_2},\theta_{j_3})dX\nonumber
\\ &&=\frac{1}{(2\pi)^3}\int_{-\pi}^\pi \int_{-\pi}^{\pi} \int_{-\pi}^{\pi}
  \bigg(  J(i\theta_1,i\theta_2;-i\theta_3)
  \nonumber\\ &&\qquad     \qquad +  J(i\theta_1,i\theta_3;-i\theta_2)+ J(i\theta_2,i\theta_3;-i\theta_1)
     +
  J(-i\theta_1,-i\theta_2;i\theta_3) \nonumber\\
&& \qquad \qquad \quad  + J(-i\theta_1,-i\theta_3;i\theta_2)
+J(-i\theta_2,-i\theta_3;i\theta_1)\nonumber \\
&&\qquad\qquad\qquad+N\big(
  J(-i\theta_1;i\theta_3)+ J(-i\theta_2;i\theta_3)+J(-i\theta_1;i\theta_2)\nonumber\\
&&\qquad \qquad \qquad \quad +
  J(-i\theta_3;i\theta_2)+
J(-i\theta_2;i\theta_1)+J(-i\theta_3;i\theta_1)\big)\nonumber
\\
&&\qquad\qquad\qquad\qquad+N^3 \bigg)
f(\theta_1,\theta_2,\theta_3) ~d\theta_1 ~d\theta_2 ~d\theta_3
\end{eqnarray}
where
\begin{eqnarray}\label{eq:IaabRMT}  &&J(\alpha_1,\alpha_2;\beta)\nonumber \\
   &&:=-e^{-\alpha_1-\alpha_2-\beta}
\int_{U(N)}\frac{\Lambda_X'}{\Lambda_X}(e^{-\alpha_1})\frac{\Lambda_X'}{\Lambda_X}(e^{-\alpha_2})
  \frac{\Lambda_{X^*}'}{\Lambda_{X^*}}(e^{-\beta})  dX \nonumber\\
&&= e^{-N(\alpha_1+\beta)}z(\alpha_1+\beta)z(-\alpha_1-\beta)
\left(\frac{z'}{z}(\alpha_2-\alpha_1)-\frac{z'}{z}(\alpha_2+\beta) \right)\nonumber\\
&&  \qquad +
e^{-N(\alpha_2+\beta)}z(\alpha_2+\beta)z(-\alpha_2-\beta)
\left(\frac{z'}{z}(\alpha_1-\alpha_2)-\frac{z'}{z}(\alpha_1+\beta)
\right)
 ,
\end{eqnarray}
and
\begin{eqnarray}
\label{eq:I1abRMT} J(\alpha;\beta):&=&e^{-\alpha-\beta}
\int_{U(N)} \frac{\Lambda_X'}{\Lambda_X}(e^{-\alpha})
  \frac{\Lambda_{X^*}'}{\Lambda_{X^*}}(e^{-\beta})  dX\nonumber\\
&=&\left(\frac{z'}{z}\right)'(\alpha+\beta)
+e^{-N(\alpha+\beta)}z(\alpha+\beta)z(-\alpha-\beta).
\end{eqnarray}
\end{theorem}

\subsection{Moments of the logarithmic derivative of the Riemann
zeta function}

Calculating correlation functions of the Riemann zeros using the
conjectural formulae for averages of ratios of zeta functions
proceeds via moments of the logarithmic derivative of the Riemann
zeta function:
\begin{eqnarray}
\label{eq:I}
&&I_r(\alpha_1,\ldots,\alpha_k;\beta_1,\ldots,\beta_{\ell})\nonumber
\\
&&\qquad\qquad:=\int_0^T
\log^r\tfrac{t}{2\pi}\frac{\zeta'}{\zeta}(\tfrac{1}{2}
+it+\alpha_1)\cdots\frac{\zeta'}{\zeta}(\tfrac{1}{2}
+it+\alpha_k)\frac{\zeta'}{\zeta}(\tfrac{1}{2} -it+\beta_1)\cdots
\frac{\zeta'}{\zeta}(\tfrac{1}{2} -it+\beta_{\ell})dt.
\end{eqnarray}
A property of these moments that will be of use to us later is
that
\begin{equation}
\label{eq:switch}
I_r(\alpha_1,\ldots,\alpha_k;\beta_1,\ldots,\beta_{\ell})=\overline{
I_r(\overline{\beta_1},\ldots,\overline{\beta_{\ell}};\overline{\alpha_1},\ldots,\overline{\alpha_k})}.
\end{equation}

For now the arguments of $I_r$ will always have positive real
parts.  In this case we have an approximate translation
invariance: if $\Re(\alpha_i),\Re(\beta_j),\Re(\alpha_i +\lambda),
\Re(\beta_j-\lambda)>0$ then
\begin{eqnarray}
\label{eq:translationinv}
&&I_r(\alpha_1+\lambda,\ldots,\alpha_k+\lambda;
\beta_1-\lambda,\ldots, \beta_{\ell}-\lambda) \\&&\qquad=
I_r(\alpha_1,\ldots,\alpha_k;\beta_1,\ldots,\beta_{\ell})
+O(|\lambda|T^{\epsilon}).\nonumber
\end{eqnarray}
This can be seen by a change of variables $t\rightarrow
t-i\lambda$ in (\ref{eq:I}) and using RH to bound
$\zeta'/\zeta(\sigma+it)$ by $t^{\epsilon}$.

In particular, to calculate the three-point correlation function,
we will need formula (\ref{eq:Iaab}) for
$I(\alpha_1,\alpha_2;\beta)$ (here we introduce the convention
that $I= I_0$).

 We will now proceed to derive that formula
using the form of the ratios conjecture \cite{kn:cfz2} for three
zeta functions over three zeta functions (with the conditions
$-\frac{1}{4} < \Re \alpha_1,\Re \alpha_2,\Re \beta< \frac{1}{4}$,
$\frac{1}{\log T} \ll\Re \gamma_1,\Re \gamma_2, \Re
\delta<\frac{1}{4}$ and  $\Im \alpha_1,\Im \alpha_2,\Im \beta,\Im
\gamma_1,\Im \gamma_2,\Im \delta\ll_{\epsilon} T^{1-\epsilon}$):
\begin{eqnarray}
&&\int_0^T
\frac{\zeta(\tfrac12+it+\alpha_1)\zeta(\tfrac12+it+\alpha_2)\zeta(
\tfrac12-it+\beta)} {\zeta(\tfrac12 +it +\gamma_1)
\zeta(\tfrac12+it +\gamma_2)\zeta(\tfrac12-it +\delta)}
dt\nonumber \\
&&\qquad \label{eq:3terms}= \int_0^T
G(\alpha_1,\alpha_2,\beta;\gamma_1,\gamma_2,\delta) +
\Big(\frac{t}{2\pi}\Big)^{-\alpha_1-\beta}
G(-\beta,\alpha_2,-\alpha_1;\gamma_1,\gamma_2,\delta) \\
&&\qquad\qquad\nonumber+
\Big(\frac{t}{2\pi}\Big)^{-\alpha_2-\beta}
G(\alpha_1,-\beta,-\alpha_2;\gamma_1,\gamma_2,\delta)\;dt +
O(T^{1/2+\epsilon}),
\end{eqnarray}
where the error term is uniform in the specified range of
parameters.  Here
\begin{eqnarray}
G(\alpha_1,\alpha_2,\beta;\gamma_1,\gamma_2,\delta)=Y(\alpha_1,\alpha_2,\beta;\gamma_1,\gamma_2,\delta)
\times A_{\zeta}(\alpha_1,\alpha_2,\beta;\gamma_1,\gamma_2,\delta)
\end{eqnarray}
and
\begin{eqnarray}
&&Y(\alpha_1,\alpha_2,\beta;\gamma_1,\gamma_2,\delta)=\frac{\zeta(1+\alpha_1+\beta)
\zeta(1+\alpha_2+\beta)
\zeta(1+\gamma_1+\delta)\zeta(1+\gamma_2+\delta)}
{\zeta(1+\alpha_1+\delta)\zeta(1+\alpha_2+\delta)
\zeta(1+\gamma_1+\beta) \zeta(1+\gamma_2+\beta)}
\end{eqnarray}
and
\begin{eqnarray}
&&A_{\zeta}(\alpha_1,\alpha_2,\beta;\gamma_1,\gamma_2,\delta)=\nonumber
\\
&&\qquad \prod_p \frac{\big(1-\frac{1}{p^{1+\alpha_1+\beta}}\big)
\big(1-\frac{1}{p^{1+\alpha_2+\beta}}\big)
\big(1-\frac{1}{p^{1+\gamma_1+\delta}}\big)
\big(1-\frac{1}{p^{1+\gamma_2+\delta}}\big)}
{\big(1-\frac{1}{p^{1+\alpha_1+\delta}}\big)
\big(1-\frac{1}{p^{1+\alpha_2+\delta}}\big)
\big(1-\frac{1}{p^{1+\gamma_1+\beta}}\big)
\big(1-\frac{1}{p^{1+\gamma_2+\beta}}\big)}\nonumber \\
&&\qquad \quad\times \sum_{{m_1+m_2+h_1+h_2=n+k}\atop
{m_1,m_2,h_1,h_2,n,k\geq 0}} \frac{ \mu(p^{h_1}) \mu(p^{h_2})
\mu(p^{k})} {p^{m_1(1/2+\alpha_1) +m_2(1/2+\alpha_2) +n(1/2+\beta)
+ h_1(1/2+\gamma_1) +h_2(1/2+\gamma_2)
+k(1/2+\delta)}}\nonumber \\
&&\qquad =\prod_p\frac{
\big(1-\frac{1}{p^{1+\gamma_1+\delta}}\big)
\big(1-\frac{1}{p^{1+\gamma_2+\delta}}\big)}
{\big(1-\frac{1}{p^{1+\alpha_1+\delta}}\big)
\big(1-\frac{1}{p^{1+\alpha_2+\delta}}\big)
\big(1-\frac{1}{p^{1+\gamma_1+\beta}}\big)
\big(1-\frac{1}{p^{1+\gamma_2+\beta}}\big)}\nonumber \\
&&\qquad \quad\times \bigg(1- p^{\beta-\delta}-\frac{1}
{p^{1+\gamma_1+\beta}} +\frac{1}{ p^{1+\gamma_1+\delta}}-
\frac{1}{ p^{1+\gamma_2+\beta}}+ \frac{1}{p^{1+\gamma_2+\delta}}  \nonumber \\
&&\qquad\qquad + \frac{1}{p^{2+\gamma_1+\gamma_2 +2\beta}}
-\frac{1}{p^{2+\gamma_1+\gamma_2+\beta+\delta}} + p^{\beta-\delta}
\big(1-\frac{1}{p^{1+\alpha_1+\beta}}\big)
\big(1-\frac{1}{p^{1+\alpha_2+\beta}}\big)\bigg).
\end{eqnarray}
The main term of (\ref{eq:3terms}) is analytic in the specified
range of parameters; the apparent poles from $Y$ cancel, as can be
checked directly, or by writing the three terms as a contour
integral as in \cite{kn:cfz2} or Section 2.5 of \cite{kn:cfkrs}.

Here $\mu(n)$ is the M{\"o}bius function and the final expression
above reflects the fact that $\mu(p^m)$ is 1 for $m=0$, it is $-1$
for $m=1$ and zero for any power of $p$ higher than the first.
Also note that
$A_{\zeta}(\alpha_1,\alpha_2,\beta;\gamma_1,\gamma_2,\beta)=1$.

We will not go through the reasoning behind this conjecture in
full, but a simpler example (one zeta function over one zeta
function) can be seen in full detail in \cite{kn:consna06} and the
original reference for the general case is \cite{kn:cfz2}.  We
will say only that the recipe for arriving at a ratios conjecture
involves replacing each zeta function in the numerator with an
approximate functional equation, those in the denominator with the
Dirichlet series for $\frac{1}{\zeta(s)}$ and then applying a
series of rules to discard all the terms in the resulting multiple
sums except for those seen above in $A_{\zeta}$.  The purpose of
$Y$ is to factor out the divergent terms in the these sums,
leaving $A_{\zeta}$ convergent for small values of
$\alpha_1,\alpha_2,$ $\beta,\gamma_1,\gamma_2$ and $\delta$.

Armed with this ratios conjecture, we want to evaluate
\begin{eqnarray}
&&I(\alpha_1,\alpha_2;\beta)\nonumber
\\
&&\qquad:=\int_0^T \frac{\zeta'}{\zeta}(\tfrac{1}{2}
+it+\alpha_1)\frac{\zeta'}{\zeta}(\tfrac{1}{2}
+it+\alpha_2)\frac{\zeta'}{\zeta}(\tfrac{1}{2}
-it+\beta)dt\nonumber \\
&&\qquad =\frac{d}{d\alpha_1} \frac{d}{d\alpha_2} \frac{d}{d\beta}
\int_0^T
\frac{\zeta(\tfrac12+it+\alpha_1)\zeta(\tfrac12+it+\alpha_2)\zeta(
\tfrac12-it+\beta)} {\zeta(\tfrac12 +it +\gamma_1)
\zeta(\tfrac12+it +\gamma_2)\zeta(\tfrac12-it +\delta)}
dt\Big|_{{{\gamma_1=\alpha_1}\atop{\gamma_2=\alpha_2}}\atop{\delta=\beta}}.
\end{eqnarray}

Examining the derivative of the first term in (\ref{eq:3terms}),
we find a great deal of cancellation upon substituting
$\gamma_1=\alpha_1$, $\gamma_2=\alpha_2$ and $\delta=\beta$, and
the only surviving terms are
\begin{eqnarray}
\label{eq:1stterm} &&\frac{d}{d\alpha_1} \frac{d}{d\alpha_2}
\frac{d}{d\beta}
G(\alpha_1,\alpha_2,\beta;\gamma_1,\gamma_2,\delta)
\Big|_{{{\gamma_1=\alpha_1}\atop{\gamma_2=\alpha_2}}\atop{\delta=\beta}}
\\
&&\qquad =\Big(\frac{\zeta'}{\zeta}\Big)'(1+\alpha_1+\beta)
\frac{d}{d\alpha_2} A_{\zeta} (\alpha_1,\alpha_2,\beta;
\alpha_1,\gamma_2 ,\beta)\Big|_{\gamma_2=\alpha_2} \nonumber \\
&&\qquad\qquad +\Big(\frac{\zeta'}{\zeta}\Big)'(1+\alpha_2+\beta)
\frac{d}{d\alpha_1} A_{\zeta} (\alpha_1,\alpha_2,\beta;
\gamma_1,\alpha_2 ,\beta)\Big|_{\gamma_1=\alpha_1} \nonumber \\
&&\qquad\qquad\qquad +\frac{d}{d\alpha_1} \frac{d}{d\alpha_2}
\frac{d}{d\beta}A_{\zeta} (\alpha_1,\alpha_2,\beta;
\gamma_1,\gamma_2
,\delta)\Big|_{{{\gamma_1=\alpha_1}\atop{\gamma_2=\alpha_2}}\atop{\delta=\beta}}\nonumber
\\
&&\qquad =\frac{d}{d\alpha_1} \frac{d}{d\alpha_2}
\frac{d}{d\beta}A_{\zeta} (\alpha_1,\alpha_2,\beta;
\gamma_1,\gamma_2
,\delta)\Big|_{{{\gamma_1=\alpha_1}\atop{\gamma_2=\alpha_2}}\atop{\delta=\beta}},\nonumber
\end{eqnarray}
because
$A_{\zeta}(\alpha_1,\alpha_2,\beta;\gamma_1,\gamma_2,\beta)=1$.

The term
\begin{eqnarray}
\label{eq:2ndterm} &&\frac{d}{d\alpha_1} \frac{d}{d\alpha_2}
\frac{d}{d\beta} \Big(\frac{t}{2\pi}\Big)^{-\alpha_1-\beta}
G(-\beta,\alpha_2,-\alpha_1;\gamma_1,\gamma_2,\delta)
\Big|_{{{\gamma_1=\alpha_1}\atop{\gamma_2=\alpha_2}}\atop{\delta=\beta}}\\
&&=\frac{d}{d\alpha_1} \frac{d}{d\alpha_2} \frac{d}{d\beta}
\Big(\frac{t}{2\pi}\Big)^{-\alpha_1-\beta}
\frac{\zeta(1-\beta-\alpha_1) \zeta(1+\alpha_2-\alpha_1)
\zeta(1+\gamma_1+\delta) \zeta(1+\gamma_2+\delta)}
{\zeta(1-\beta+\delta)\zeta(1+\alpha_2+\delta)
\zeta(1+\gamma_1-\alpha_1) \zeta(1+\gamma_2-\alpha_1)}\nonumber \\
&&\qquad\qquad \times
A_{\zeta}(-\beta,\alpha_2,-\alpha_1;\gamma_1,\gamma_2,\delta)
 \Big|_{{{\gamma_1=\alpha_1}\atop{\gamma_2=\alpha_2}}\atop{\delta=\beta}}
 \nonumber
 \end{eqnarray}
 quickly reduces to
\begin{eqnarray}
&&\frac{d}{d\alpha_2} \Big(\frac{t}{2\pi}\Big)^{-\alpha_1-\beta}
\frac{\zeta(1-\beta-\alpha_1) \zeta(1+\alpha_2-\alpha_1)
\zeta(1+\alpha_1+\beta) \zeta(1+\gamma_2+\beta)}
{\zeta(1+\alpha_2+\beta) \zeta(1+\gamma_2-\alpha_1)}\nonumber\\
&&\qquad\qquad\times
A_{\zeta}(-\beta,\alpha_2,-\alpha_1;\alpha_1,\gamma_2,\beta)
 \Big|_{\gamma_2=\alpha_2}
\end{eqnarray}
 because we see immediately that the factors
$\frac{1}{\zeta(1-\beta+\delta)}$ and
$\frac{1}{\zeta(1+\gamma_1-\alpha_1)}$ cause the entire term to
evaluate as zero upon the substitution $\gamma_1=\alpha_1$ and
$\delta=\beta$ unless the $\beta$ and $\alpha_1$ derivatives are
performed on these factors.  The final differentiation with
respect to $\alpha_2$ shows us that (\ref{eq:2ndterm}) equals
\begin{eqnarray}
&&\Big(\frac{t}{2\pi}\Big)^{-\alpha_1-\beta}
\zeta(1-\beta-\alpha_1)
\zeta(1+\beta+\alpha_1) \nonumber \\
&&\quad\times\bigg(\Big( \frac{\zeta'}{\zeta}
(1+\alpha_2-\alpha_1) - \frac{\zeta'}{\zeta}(1+\alpha_2
+\beta)\Big) A_{\zeta}(-\beta,\alpha_2,-\alpha_1;
\alpha_1,\alpha_2,\beta)\nonumber \\
&&\qquad\qquad +\frac{d}{d\alpha_2} A_{\zeta} (-\beta,
\alpha_2,-\alpha_1;\alpha_1,\gamma_2,\beta)\Big|_{\gamma_2=\alpha_2}\bigg).
\end{eqnarray}

Similarly,
\begin{eqnarray} \label{eq:3rdterm}&&\frac{d}{d\alpha_1} \frac{d}{d\alpha_2} \frac{d}{d\beta}
\Big(\frac{t}{2\pi}\Big)^{-\alpha_2-\beta}
G(\alpha_1,-\beta,-\alpha_2;\gamma_1,\gamma_2,\delta)
\Big|_{{{\gamma_1=\alpha_1}\atop{\gamma_2=\alpha_2}}\atop{\delta=\beta}}\\
&&=\Big(\frac{t}{2\pi}\Big)^{-\alpha_2-\beta}
\zeta(1-\beta-\alpha_2) \zeta(1+\beta+\alpha_2) \nonumber
\\
&&\quad \times\bigg(\Big( \frac{\zeta'}{\zeta}
(1+\alpha_1-\alpha_2) - \frac{\zeta'}{\zeta}(1+\alpha_1
+\beta)\Big) A_{\zeta}(\alpha_1,-\beta,-\alpha_2;
\alpha_1,\alpha_2,\beta) \nonumber \\
&&\qquad\qquad +\frac{d}{d\alpha_1} A_{\zeta}
(\alpha_1,-\beta,-\alpha_2;\gamma_1,\alpha_2,\beta)\Big|_{\gamma_1=\alpha_1}\bigg).\nonumber
\end{eqnarray}

Finally, some manipulation of the prime product $A_\zeta$ (for
which Mathematica is very helpful) shows us (where $A$, $P$ and
$Q$ refer to equations (\ref{eq:Aprime}), (\ref{eq:Pprime}) and
(\ref{eq:Qprime}), respectively) that
\begin{eqnarray}
&&A_{\zeta}(-\beta,\alpha_2,-\alpha_1;
\alpha_1,\alpha_2,\beta)=A(\alpha_1+\beta)\\
&&\quad =\prod_p
\frac{(1-\tfrac{1}{p^{1+\alpha_1+\beta}})(1-\tfrac{2}{p}+\tfrac{1}{p^{1+\alpha_1+\beta}})}
{(1-\tfrac{1}{p})^2},\nonumber \\&&\nonumber \\
&&\frac{d}{d\alpha_2} A_{\zeta} (-\beta,
\alpha_2,-\alpha_1;\alpha_1,\gamma_2,\beta)\Big|_{\gamma_2=\alpha_2}=P(\alpha_1+\beta,\alpha_2+\beta)\\
&&\quad =
\prod_p\frac{\Big(1-\frac{1}{p^{1+\alpha_1+\beta}}\Big)\Big(1-\frac{2}{p}+\frac{1}{p^{1+\alpha_1+\beta}}\Big)
} {\Big(1-\frac{1}{p}\Big)^2}\nonumber \\
&&\qquad\times
\sum_p\frac{\Big(1-\frac{1}{p^{\alpha_1+\beta}}\Big)
\Big(1-\frac{1}{p^{\alpha_1+\beta}}-\frac{1}{p^{\alpha_2+\beta}}
+\frac{1}{p^{1+\alpha_2+\beta}}\Big) \log p}
{\Big(\frac{1}{p^{1-\alpha_1+\alpha_2}}-1\Big)
\Big(1-\frac{1}{p^{1+\alpha_2+\beta}}\Big) \Big(
1-\frac{2}{p}+\frac{1}{p^{1+\alpha_1+\beta}}\Big)p^{2-\alpha_1+\alpha_2}},\nonumber
\\&&\nonumber \\
&&\frac{d}{d\alpha_1} A_{\zeta}
(\alpha_1,-\beta,-\alpha_2;\gamma_1,\alpha_2,\beta)\Big|_{\gamma_1=\alpha_1}=P(\alpha_2+\beta,\alpha_1+\beta),\\
&&\nonumber \\
&&\frac{d}{d\alpha_1} \frac{d}{d\alpha_2}
\frac{d}{d\beta}A_{\zeta} (\alpha_1,\alpha_2,\beta;
\gamma_1,\gamma_2 ,\delta)\Big|_{{{\gamma_1=\alpha_1}
\atop{\gamma_2=\alpha_2}}\atop{\delta=\beta}}=Q(\alpha_1+\beta,\alpha_2+\beta)\\
&&\quad =-\sum_p \frac{\log^3 p}
{p^{2+\alpha_1+\alpha_2+2\beta}(1-\tfrac{1}{p^{1+\alpha_1+\beta}})(
1-\tfrac{1}{p^{1+\alpha_2+\beta}})}.\nonumber
\end{eqnarray}
Substituting these expressions into (\ref{eq:1stterm}),
(\ref{eq:2ndterm}) and (\ref{eq:3rdterm}), we arrive at
(\ref{eq:Iaab}) as expected.

The other version of (\ref{eq:I}) that we need is
$I_1(\alpha;\beta)$, as given in (\ref{eq:I1ab}).  This
calculation has in fact already been carried out in
\cite{kn:consna06} using the two zeta functions over two zeta
functions ratios conjecture in a manner very similar to the
three-over-three calculation above.  That is
\cite{kn:cfz1,kn:cfz2}, with $-\frac{1}{4} < \Re \alpha, \Re
\beta< \frac{1}{4}$, $\frac{1}{\log T} \ll \Re \gamma, \Re
\delta<\frac{1}{4}$ and $\Im \alpha, \Im \beta,\Im \gamma, \Im
\delta\ll_{\epsilon} T^{1-\epsilon}$ (for every $\epsilon>0$), the
ratios conjecture states
\begin{eqnarray}
&&\int_0^T
 \nonumber\frac{\zeta(s+\alpha)\zeta(1-s+\beta)}{\zeta(s+\gamma)\zeta(1-s+\delta)}~dt\\
 &&\qquad=\int_0^T \left(\frac{\zeta(1+\alpha+\beta)\zeta(1+\gamma+\delta)}
{\zeta(1+\alpha+\delta)\zeta(1+\beta+\gamma)}A_\zeta(\alpha, \beta; \gamma, \delta) \right. \\
&& \qquad\quad \left. +\left(\frac {t}{2
\pi}\right)^{-\alpha-\beta}
\frac{\zeta(1-\alpha-\beta)\zeta(1+\gamma+\delta)}
{\zeta(1-\beta+\delta)\zeta(1-\alpha+\gamma)}A_\zeta(-\beta,
-\alpha;  \gamma, \delta) \right)~dt
+O\left(T^{1/2+\epsilon}\right),\nonumber
\end{eqnarray}
where
\begin{eqnarray}
A_\zeta(\alpha, \beta; \gamma, \delta)=\prod_p
\frac{\left(1-\frac{1}{p^{1+\gamma+\delta}}\right)
\left(1-\frac{1}{p^{1+\beta+\gamma}}
-\frac{1}{p^{1+\alpha+\delta}}+\frac{1}{p^{1+\gamma+\delta}}\right)}
{\left(1-\frac{1}{p^{1+\beta+\gamma}}\right)
\left(1-\frac{1}{p^{1+\alpha+\delta}}\right)}.
\end{eqnarray}
This implies that
\begin{eqnarray}
&& \int_0^T
\frac{\zeta'}{\zeta}(s+\alpha)\frac{\zeta'}{\zeta}(1-s+\beta)~dt
=\int_0^T\left( \left(\frac{\zeta'}{\zeta}\right)'(1+\alpha+\beta)+\right. \\
&&\quad \left(\frac{t}{2\pi}\right)^ {-\alpha-\beta}
\zeta(1+\alpha+\beta)\zeta(1-\alpha-\beta)
\prod_p\frac{(1-\frac{1}{p^{1+\alpha+\beta}})
(1-\frac 2 p +\frac{1}{p^{1+\alpha+\beta}})}{(1-\frac 1 p )^2}\nonumber\\
&& \quad \left. -\sum_p \left(\frac{\log
p}{(p^{1+\alpha+\beta}-1)} \right)^2\right)~dt
+O(T^{1/2+\epsilon}),\nonumber
\end{eqnarray}
provided that $\frac{1}{\log T}\ll\Re \alpha, \Re \beta
<\frac{1}{4}$.  The ratios conjecture recipe can incorporate the
$\log \tfrac{t}{2\pi}$ factor in $I_1(\alpha;\beta)$ without any
alteration, giving (\ref{eq:I1ab}).

\subsection{Triple correlation as a contour integral}
\label{sect:contour}

We start with the a triple sum over zeros of the Riemann zeta
function, where the zeros do not have to be distinct.  Using
Cauchy's residue theorem, the triple sum can be written as a
triple contour integral where each contour is a rectangle
enclosing the zeros $1/2+i\gamma$ (assuming the Riemann
Hypothesis) with heights $0\leq \gamma<T$:
\begin{eqnarray}
&&\sum_{0\leq \gamma_1, \gamma_2,\gamma_3<T}
g(\gamma_1,\gamma_2,\gamma_3)\\
&&\qquad\qquad=\frac{1}{(2\pi i)^3}\oint\oint\oint
g(-i(x-1/2),-i(y-1/2),-i(z-1/2)) \frac{\zeta'}{\zeta}(x)
\frac{\zeta'}{\zeta}(y)\frac{\zeta'}{\zeta}(z) dx dy dz.\nonumber
\end{eqnarray}

In investigating the 3-point correlation of the Riemann zeros, we
are interested in the relative spacing between triples of zeros,
so we can assume that the test function is translation invariant.
Thus we define a function $f$ satisfying the conditions
\begin{eqnarray} \label{eq:fconditions}
&&f(x,y) \text{ is holomorphic for }  |\Im x| <2, |\Im y| <2\\
 && \text{and that } f(x,y)\ll
1/(1+|x|^2+|y|^2) \text{ as } |x| \text{ or } |y| \to
\infty.\nonumber
\end{eqnarray}
  Thus we have
\begin{eqnarray}
\label{eq:triplecontour} \sum_{0\leq \gamma_1,
\gamma_2,\gamma_3<T}
f(\gamma_1-\gamma_2,\gamma_1-\gamma_3)&=&\frac{1}{(2\pi i)^3}\oint
\oint\oint  \frac{\zeta'}{\zeta}(z_1)
\frac{\zeta'}{\zeta}(z_2)\frac{\zeta'}{\zeta}(z_3)\nonumber \\
&&\quad\times f\big(-i(z_1-z_2),-i(z_1-z_3)\big) dz_1 dz_2dz_3,
\end{eqnarray}
where the contours are rectangles with corners at the points
$(a,0)$, $(a,iT)$, $(1-b,iT)$ and $(1-b,0)$ with $1/2<a,b<1$.  We
distinguish $a$ and $b$ only for ease of following the
manipulations in the following calculations.

The horizontal portions of the contour of integration can be
chosen so that the integral along them is negligible, so we
concentrate on the vertical sides of the contours.  This makes
(\ref{eq:triple}) the sum of eight integrals $J_1,\ldots,J_8$
which will be met one by one below (a subscript $a$ on an integral
indicates integration from $(a,0)$ to $(a,iT)$).

First we have
\begin{eqnarray}
\label{eq:J1} J_1&=&\frac{1}{(2\pi i)^3} \int_a\int_a\int_a
\frac{\zeta'}{\zeta}(z_1)
\frac{\zeta'}{\zeta}(z_2)\frac{\zeta'}{\zeta}(z_3)\nonumber\\
&&\qquad\times
f\big(-i(z_1-z_2),-i(z_1-z_3)\big) dz_1 dz_2dz_3\nonumber \\
&=&O(T^{\epsilon}).
\end{eqnarray}
The final line is true because all three contours can be moved to
the right (assuming the Riemann Hypothesis) where
$\frac{\zeta'}{\zeta}$ converges and can be integrated term by
term (the pole at 1 doesn't contribute more than a constant).

Next we examine
\begin{eqnarray}
\label{eq:J2} J_2&:=&-\frac{1}{(2\pi i)^3} \int_a\int_a\int_{1-b}
\frac{\zeta'}{\zeta}(z_1)
\frac{\zeta'}{\zeta}(z_2)\frac{\zeta'}{\zeta}(z_3)\nonumber\\
&&\qquad\times f\big(-i(z_1-z_2),-i(z_1-z_3)\big) dz_1 dz_2dz_3.
\end{eqnarray}
We use the functional equation
\begin{equation}
\label{eq:logderivfe}\frac{\zeta'}{\zeta}(s)=\frac{\chi'}{\chi}(s)-\frac{\zeta'}{\zeta}(1-s)
\end{equation}
and obtain
\begin{eqnarray}
\label{eq:J2a} J_2&=&-\frac{1}{(2\pi i)^3} \int_a\int_a\int_{1-b}
\Big(\frac{\chi'}{\chi}(z_1)- \frac{\zeta'}{\zeta}(1-z_1)\Big)
\frac{\zeta'}{\zeta}(z_2)\frac{\zeta'}{\zeta}(z_3)\nonumber\\
&&\qquad\times f\big(-i(z_1-z_2),-i(z_1-z_3)\big) dz_1 dz_2dz_3.
\end{eqnarray}
For the term with $\frac{\chi'}{\chi}(z_1)$ we can shift the
integrals to the right, as we did with $J_1$, and so see that the
contribution is only $O(T^{\epsilon})$.  A similar manoeuver
cannot be done with the integral containing the three
$\frac{\zeta'}{\zeta}$, however, because the $z_1$ integral is to
the left of the critical line.  So,
\begin{eqnarray}
\label{eq:J2b} J_2&=&\frac{1}{(2\pi)^3} \int_0^T\int_0^T\int_0^T
 \frac{\zeta'}{\zeta}(b-it_1)
\frac{\zeta'}{\zeta}(a+it_2)\frac{\zeta'}{\zeta}(a+it_3)\nonumber\\
&&\qquad\times f\big(t_1-t_2-i(1-b-a),t_1-t_3-i(1-b-a)\big) dt_1
dt_2dt_3+O(T^{\epsilon}).
\end{eqnarray}
Now, letting $v_1=t_2-t_1$
and $v_2=t_3-t_1$, we have
\begin{eqnarray}
\label{eq:J2b1}
J_2&=&\frac{1}{(2\pi)^3}
\int_0^T\int_{-t_1}^{T-t_1}\int_{-t_1}^{T-t_1}
 \frac{\zeta'}{\zeta}(b-it_1)
\frac{\zeta'}{\zeta}(a+i(v_1+t_1))\frac{\zeta'}{\zeta}(a+i(v_2+t_1))\nonumber\\
&&\qquad\times f\big(-v_1-i(1-b-a),-v_2-i(1-b-a)\big) dv_1
dv_2dt_1+O(T^{\epsilon})\nonumber \\
&=&\frac{1}{(2\pi)^3}
\int_{-T}^T\int_{-T}^{T}f\big(-v_1-i(1-b-a),-v_2-i(1-b-a)\big)\nonumber\\
&&\qquad\times \int_{0}^{T}
 \frac{\zeta'}{\zeta}(b-it_1)
\frac{\zeta'}{\zeta}(a+i(v_1+t_1))\frac{\zeta'}{\zeta}(a+i(v_2+t_1))
dt_1 dv_1 dv_2+O(T^{\epsilon}).
\end{eqnarray}
In this last line we have switched the order of integration of
$t_1$ with $v_1$ and $v_2$.  The range of the innermost integral
should really be from $\max\{0,-v_1,-v_2\}$ to
$\min\{T,T-v_1,T-v_2\}$.  However, since we are assuming that
$f(x,y)$ decays fast, see (\ref{eq:fconditions}), it is not hard
to show that extending the interval to $(0,T)$ as we have done
above, and will do again in the following integrals $J_3$ to
$J_8$, incurs only an error that is a power of $\log T$ and so can
be incorporated into the error term $O(T^{\epsilon})$.

Now all that is left is to tidy up the expression for $J_2$, so we
shift the contours of integration off the real axis and apply the
definition of $I(\alpha_1,\alpha_2;\beta)$ from (\ref{eq:I}) with
$r=0$:
\begin{eqnarray}
\label{eq:J2c}
 J_2&=&\frac{1}{(2\pi)^3}
\int_{-T-i(1-b-a)}^{T-i(1-b-a)}\int_{-T-i(1-b-a)}^{T-i(1-b-a)}f(v_1,v_2)\nonumber\\
&&\qquad\times \int_{0}^{T}
 \frac{\zeta'}{\zeta}(b-it_1)
\frac{\zeta'}{\zeta}(1+it_1-b-iv_1)\frac{\zeta'}{\zeta}(1+it_1-b-iv_2)
dt_1 dv_1 dv_2+O(T^{\epsilon})\nonumber \\
&=&\frac{1}{(2\pi)^3}
\int_{-T-i(1-b-a)}^{T-i(1-b-a)}\int_{-T-i(1-b-a)}^{T-i(1-b-a)}f(v_1,v_2)I(-iv_1,-iv_2;0)dv_1dv_2+
O(T^{\epsilon}).
\end{eqnarray}
To simplify the last line we have used (\ref{eq:translationinv}).

We use exactly the same sequence of manipulations to obtain
similar expressions for $J_3$ and $J_4$:
\begin{eqnarray}
\label{eq:J3} J_3&:=&-\frac{1}{(2\pi i)^3}
\int_a\int_{1-b}\int_{a} \frac{\zeta'}{\zeta}(z_1)
\frac{\zeta'}{\zeta}(z_2)\frac{\zeta'}{\zeta}(z_3)\nonumber\\
&&\qquad\times f\big(-i(z_1-z_2),-i(z_1-z_3)\big) dz_1 dz_2dz_3
\nonumber
\\
&=& \frac{1}{(2\pi)^3}
\int_{-T}^{T}\int_{-T+i(1-b-a)}^{T+i(1-b-a)}f(v_1,v_2)I(0,-iv_2;iv_1)dv_1dv_2
+O(T^{\epsilon})
\end{eqnarray}
and
\begin{eqnarray}
\label{eq:J4} J_4&:=&-\frac{1}{(2\pi i)^3}
\int_{1-b}\int_{a}\int_{a} \frac{\zeta'}{\zeta}(z_1)
\frac{\zeta'}{\zeta}(z_2)\frac{\zeta'}{\zeta}(z_3)\nonumber\\
&&\qquad\times f\big(-i(z_1-z_2),-i(z_1-z_3)\big) dz_1 dz_2dz_3
\nonumber
\\
&=& \frac{1}{(2\pi)^3}
\int_{-T+i(1-b-a)}^{T+i(1-b-a)}\int_{-T}^{T}f(v_1,v_2)I(0,-iv_1;iv_2)dv_1dv_2
+O(T^{\epsilon}).
\end{eqnarray}

The integral $J_5$ throws up something slightly different.  We
start off in an identical manner, replacing
$\frac{\zeta'}{\zeta}(z_2)$ and $\frac{\zeta'}{\zeta}(z_3)$ by
their functional equation (\ref{eq:logderivfe}).  This results in
four terms, one of which contributes no more than
$O(T^{\epsilon})$.  We are left with three terms:
\begin{eqnarray}
\label{eq:J5} J_5&:=&\frac{1}{(2\pi i)^3}
\int_{1-b}\int_{1-b}\int_{a} \frac{\zeta'}{\zeta}(z_1)
\frac{\zeta'}{\zeta}(z_2)\frac{\zeta'}{\zeta}(z_3)\nonumber\\
&&\qquad\times f\big(-i(z_1-z_2),-i(z_1-z_3)\big) dz_1 dz_2dz_3
\nonumber
\\
&=& \frac{1}{(2\pi)^3}\int_0^T\int_0^T\int_0^T
\Bigg(\frac{\zeta'}{\zeta}(a+it_1)
\frac{\zeta'}{\zeta}(b-it_2)\frac{\zeta'}{\zeta}(b-it_3)\nonumber\\
&&\qquad\qquad -
\frac{\chi'}{\chi}(1-b+it_2)\frac{\zeta'}{\zeta}(a+it_1)\frac{\zeta'}{\zeta}(b-it_3)\nonumber
\\
&&\qquad\qquad\qquad
-\frac{\chi'}{\chi}(1-b+it_3)\frac{\zeta'}{\zeta}(a+it_1)\frac{\zeta'}{\zeta}(b-it_2)\Bigg)\nonumber
\\
&&\qquad \times f\big(t_1-t_2+i(1-a-b),t_1-t_3+i(1-a-b)\big)
dt_1dt_2dt_3 +O(T^{\epsilon})
\nonumber \\
&=&\frac{1}{(2\pi)^3}\int_0^T\int_{-t_1}^{T-t_1}\int_{-t_1}^{T-t_1}
\Bigg(\frac{\zeta'}{\zeta}(a+it_1)
\frac{\zeta'}{\zeta}(b-i(v_1+t_1))\frac{\zeta'}{\zeta}(b-i(v_2+t_1))\nonumber\\
&&\qquad\qquad -
\frac{\chi'}{\chi}(1-b+i(v_1+t_1))\frac{\zeta'}{\zeta}(a+it_1)\frac{\zeta'}{\zeta}(b-i(v_2+t_1))\nonumber
\\
&&\qquad\qquad\qquad
-\frac{\chi'}{\chi}(1-b+i(v_2+t_1))\frac{\zeta'}{\zeta}(a+it_1)\frac{\zeta'}{\zeta}(b-i(v_1+t_1))\Bigg)\nonumber
\\
&&\qquad \times f\big(-v_1+i(1-a-b),-v_2+i(1-a-b)\big)
dv_1dv_2dt_1 +O(T^{\epsilon}),
\end{eqnarray}
 where in the last line we have made our usual change of
variables $v_1=t_2-t_1$ and $v_2=t_3-t_1$.

We now note the asymptotic for $\frac{\chi'}{\chi}$:
\begin{eqnarray}
\label{eq:chichi}
\frac{\chi'}{\chi}(1/2+it)=-\log\frac{|t|}{2\pi}\left(1+O\left(\frac
1 {|t|}\right)\right).
\end{eqnarray}
Since $f$ is very small when $v_1$ or $v_2$ are large, we replace
$\frac{\chi'}{\chi}(1-b+i(v_j+z_1))$ with $-\log
\tfrac{z_1}{2\pi}$ and obtain
\begin{eqnarray}
\label{eq:J5a} J_5&=&\frac{1}{(2\pi)^3}\int_{-T}^{T}\int_{-T}^{T}
f\big(-v_1+i(1-a-b),-v_2+i(1-a-b)\big) \nonumber \\
&&\qquad \times \int_0^T \Bigg(\frac{\zeta'}{\zeta}(a+it_1)
\frac{\zeta'}{\zeta}(b-i(v_1+t_1))\frac{\zeta'}{\zeta}(b-i(v_2+t_1))\nonumber\\
&&\qquad\qquad
+\log\tfrac{t_1}{2\pi}\frac{\zeta'}{\zeta}(a+it_1)\frac{\zeta'}{\zeta}(b-i(v_2+t_1))\nonumber
\\
&&\qquad\qquad\qquad +\log
\tfrac{t_1}{2\pi}\frac{\zeta'}{\zeta}(a+it_1)\frac{\zeta'}{\zeta}(b-i(v_1+t_1))\Bigg)\nonumber
dt_1 dv_1dv_2+O(T^{\epsilon}),
\end{eqnarray}
where we have extended the range of $t_1$, after exchanging the
order of integration, by the same argument as for $J_2$.

Now we can write $J_5$ as a double integral along contours running
just below the real axis, in a similar form to $J_2$, $J_3$ and
$J_4$ above, but we introduce $I_1(\alpha;\beta)$ (see
(\ref{eq:I}) and (\ref{eq:I1ab})) to arrive at
\begin{eqnarray}
J_5&=& \frac{1}{(2\pi)^3}
\int_{-T+i(1-a-b)}^{T+i(1-a-b)}\int_{-T+i(1-a-b)}^{T+i(1-a-b)}
f(v_1,v_2) \Big( I(0;iv_1,iv_2)\nonumber
\\
&& \qquad\qquad\qquad+I_1(0;iv_2)+I_1(0;iv_1)\Big)
dv_1dv_2+O(T^{\epsilon}).
\end{eqnarray}

Proceeding exactly as for $J_5$, we obtain similar expressions for
$J_6$ and $J_7$:
\begin{eqnarray}
\label{eq:J6} J_6&:=&\frac{1}{(2\pi i)^3}
\int_{1-b}\int_{a}\int_{1-b} \frac{\zeta'}{\zeta}(z_1)
\frac{\zeta'}{\zeta}(z_2)\frac{\zeta'}{\zeta}(z_3)\nonumber\\
&&\qquad\times f\big(-i(z_1-z_2),-i(z_1-z_3)\big) dz_1 dz_2dz_3
\nonumber
\\
&=& \frac{1}{(2\pi)^3}
\int_{-T}^{T}\int_{-T-i(1-a-b)}^{T-i(1-a-b)} f(v_1,v_2) \Big(
I(-iv_1;0,iv_2)\nonumber
\\
&& \qquad\qquad\qquad+I_1(-iv_1;iv_2)+I_1(-iv_1;0)\Big)
dv_1dv_2+O(T^{\epsilon})
\end{eqnarray}
and
\begin{eqnarray}
\label{eq:J7} J_7&:=&\frac{1}{(2\pi i)^3}
\int_{a}\int_{1-b}\int_{1-b} \frac{\zeta'}{\zeta}(z_1)
\frac{\zeta'}{\zeta}(z_2)\frac{\zeta'}{\zeta}(z_3)\nonumber\\
&&\qquad\times f\big(-i(z_1-z_2),-i(z_1-z_3)\big) dz_1 dz_2dz_3
\nonumber
\\
&=& \frac{1}{(2\pi)^3}\int_{-T-i(1-a-b)}^{T-i(1-a-b)}
\int_{-T}^{T} f(v_1,v_2) \Big( I(-iv_2;0,iv_1)\nonumber
\\
&& \qquad\qquad\qquad+I_1(-iv_2;iv_1)+I_1(-iv_2;0)\Big)
dv_1dv_2+O(T^{\epsilon}).
\end{eqnarray}

We are just left with the integral
\begin{eqnarray}\label{eq:J8}
J_8&:=&-\frac{1}{(2\pi i)^3} \int_{1-b}\int_{1-b}\int_{1-b}
\frac{\zeta'}{\zeta}(z_1)
\frac{\zeta'}{\zeta}(z_2)\frac{\zeta'}{\zeta}(z_3)\nonumber\\
&&\qquad\times f\big(-i(z_1-z_2),-i(z_1-z_3)\big) dz_1 dz_2dz_3
\end{eqnarray}
to evaluate.  Once each $\frac{\zeta'}{\zeta}(z)$ is replaced by
its functional equation (\ref{eq:logderivfe}), any term with at
least one $\frac{\zeta'}{\zeta}(1-z)$ in it can be shown to be
size just $O(T^{\epsilon})$ by shifting the contour far to the
left.  This leaves us with
\begin{eqnarray}
\label{eq:J8a} J_8&=&-\frac{1}{(2\pi)^3} \int_0^T \int_0^T
\int_0^T \frac{\chi'}{\chi}(\tfrac{1}{2} +it_1)
\frac{\chi'}{\chi}(\tfrac{1}{2} +it_2)
\frac{\chi'}{\chi}(\tfrac{1}{2} +it_3)\nonumber
\\
&& \qquad\qquad\qquad\times f(t_1-t_2,t_1-t_3) dt_1 dt_2 dt_3
+O(T^{\epsilon}).
\end{eqnarray}
Here the contours have all been moved onto the half-line as there
are no longer zeros of zeta functions to avoid. We replace
$t_2-t_1$ with $v_1$ and replace $t_3-t_1$ with $v_2$ and
substitute $-\log \tfrac{t_1}{2\pi}$ for each $\frac{\chi'}{\chi}$
factor with the help of (\ref{eq:chichi}), as we did in the
discussion of $J_5$, and so obtain
\begin{eqnarray}
\label{eq:J8b} J_8&=& \frac{1}{(2\pi)^3} \int_{-T}^T \int_{-T}^T
f(v_1,v_2) \int_0^T \log^3 \tfrac{t}{2\pi} dt \; dv_1
dv_2+O(T^{\epsilon}).
\end{eqnarray}

The sum of integrals $J_1,\ldots,J_8$ gives us (see
(\ref{eq:triplePV})) the expression for
\begin{equation}
\sum_{0\leq \gamma_1, \gamma_2,\gamma_3<T}
f(\gamma_1-\gamma_2,\gamma_1-\gamma_3).
\end{equation}
Such a triple sum over zeros necessarily contains terms where two
zeros are identical.  These are essentially two-point
correlations, rather than three-point statistics.  Similarly,
terms where all three zeros are identical are just one-point
statistics.  To remove these lower-order correlations, one looks
at the sum
\begin{equation}
\sum_{0<\gamma_1\neq\gamma_2\neq\gamma_3<T} f(\gamma_1-\gamma_2,
\gamma_1-\gamma_3)
\end{equation}
instead.  As will be shown in the next section, simply rewriting
as an integral on the real line each of the integrals above that
occurs with a contour on $(-T-i(1-a-b),T-i(1-a-b))$ or
$(-T+i(1-a-b),T+i(1-a-b))$ gives us the purely triple correlation,
and so we obtain the final result (\ref{eq:triple}).  The
contributions from the poles that we meet as the contours are
shifted to the real axes are the source of the two- and one-point
correlation terms.

\subsection{Contributions from lower-order correlations}

Collecting together the results from the previous section we have
\begin{eqnarray}
\label{eq:triplePV} &&\sum_{0\leq \gamma_1, \gamma_2,\gamma_3<T}
f(\gamma_1-\gamma_2,\gamma_1-\gamma_3)= \frac{1}{(2\pi)^3}
\Bigg(\int_{-T}^T \int_{-T}^T f(v_1,v_2) \int_0^T \log^3
\tfrac{t}{2\pi} dt \; dv_1 dv_2\nonumber \\
&&\qquad
+\int_{-T+i(1-a-b)}^{T+i(1-a-b)}\int_{-T+i(1-a-b)}^{T+i(1-a-b)}
f(v_1,v_2) \Big( I(0;iv_1,iv_2)\nonumber
\\
&& \qquad\qquad\qquad+I_1(0;iv_2)+I_1(0;iv_1)\Big)
dv_1dv_2\nonumber \\
&&\qquad +\int_{-T}^{T}\int_{-T-i(1-a-b)}^{T-i(1-a-b)} f(v_1,v_2)
\Big( I(-iv_1;0,iv_2)\nonumber
\\
&& \qquad\qquad\qquad+I_1(-iv_1;iv_2)+I_1(-iv_1;0)\Big) dv_1dv_2
\nonumber \\
&&\qquad +\int_{-T-i(1-a-b)}^{T-i(1-a-b)} \int_{-T}^{T} f(v_1,v_2)
\Big( I(-iv_2;0,iv_1)\nonumber
\\
&& \qquad\qquad\qquad+I_1(-iv_2;iv_1)+I_1(-iv_2;0)\Big)
dv_1dv_2\nonumber \\
&&\qquad
+\int_{-T-i(1-b-a)}^{T-i(1-b-a)}\int_{-T-i(1-b-a)}^{T-i(1-b-a)}f(v_1,v_2)I(-iv_1,-iv_2;0)dv_1dv_2\nonumber
\\
&&\qquad+
\int_{-T}^{T}\int_{-T+i(1-b-a)}^{T+i(1-b-a)}f(v_1,v_2)I(0,-iv_2;iv_1)dv_1dv_2\nonumber
\\
&&\qquad +
\int_{-T+i(1-b-a)}^{T+i(1-b-a)}\int_{-T}^{T}f(v_1,v_2)I(0,-iv_1;iv_2)dv_1dv_2\Bigg)
+O(T^{\epsilon}).
\end{eqnarray}

The goal of this section is to move all the contours of
integration above onto the real axis, evaluating the contributions
from the poles encountered during this process.  The resulting
integrals along the real axis will then be computable, for the
purposes of Figure ???? for example, as principal value integrals.
Elegantly, the terms resulting from the residue at the various
poles will yield the contribution to (\ref{eq:triplePV}) from
lower-order correlations between Riemann zeros.

Consider, as an example of the method, the term (coming from
$J_5$)
\begin{eqnarray}
\label{eq:sq1}
\int_{-T+i(1-a-b)}^{T+i(1-a-b)}\int_{-T+i(1-a-b)}^{T+i(1-a-b)}
f(v_1,v_2)  I(iv_1,iv_2;0) dv_1dv_2
\end{eqnarray}
from (\ref{eq:triplePV}).  (Here we have switched the order of the
arguments of $I$ using (\ref{eq:switch}).)  Treating the integrand
first as a function of $v_2$ alone (leaving $v_1$ fixed and
non-zero), we see that the $v_2$ contour lies below the real axis
(since $a>1/2$ and $b>1/2$) and there is a pole at $v_2=0$, so
shifting the contour onto the axis results in a principal value
integral in $v_2$ plus $i\pi$ times the residue of that pole.
Using the definition of the sums and products over primes found in
(\ref{eq:Aprime}) to (\ref{eq:Pprime}), expanding around $v_2=0$
gives us
\begin{eqnarray}
\label{eq:primeexpand} Q(iv_1,iv_2)&=&O(1)\\
A(iv_2)&=& 1+O(v_2^2)\\
P(iv_2,iv_1)&=&iv_2\sum_p \frac{\log^2
p}{(p^{1+iv_1}-1)^2}+O(v_2^2)\\
P(iv_1,iv_2)&=&O(1).
\end{eqnarray}
These, in conjunction with the poles of $I(iv_1,iv_2;0)$ at
$v_2=0$ resulting from $\zeta(1-iv_2)\zeta(1+iv_2)$ and
$\frac{\zeta'}{\zeta}(1+iv_2)$, allows us to express
(\ref{eq:sq1}) as
\begin{eqnarray}
\label{eq:sq1a} &&\int_{-T}^T\int_{-T+i(1-a-b)}^{T+i(1-a-b)}
f(v_1,v_2)I(iv_1,iv_2;0)dv_1dv_2 \nonumber \\
&&\qquad\qquad+\pi\int_{-T+i(1-a-b)}^{T+i(1-a-b)}\int_0^T \bigg(
\Big(\frac{\zeta'}{\zeta}\Big)'(1+iv_1) +
\big(\frac{t}{2\pi}\big)^{-iv_1}
\zeta(1+iv_1)\zeta(1-iv_1)A(iv_1) \nonumber \\
&&\qquad\qquad\qquad-\sum_p \frac{\log^2 p} {(p^{1+iv_1}-1)^2}
\bigg) f(v_1,0)dt \; dv_1,
\end{eqnarray}
where the $v_2$ integral is understood as a principal value
integral.

 Now we move the $v_1$ integral onto the real axis. In
this case also we encounter a pole at $v_1=0$. The first integral
in (\ref{eq:sq1a}) yields a residue contribution essentially
identical to that from the $v_2$ pole treated above, and the
second integral has a pole at $v_1=0$ with residue $-i\log
\tfrac{t}{2\pi}$.  We take $i \pi$ times the contributions from
these poles and the end result is
\begin{eqnarray}
\label{eq:sq1b} &&\int_{-T}^{T} \int_{-T}^T
f(v_1,v_2)I(iv_1,iv_2;0)dv_1dv_2 \nonumber \\
&&\qquad\qquad+\pi\int_{-T}^T f(v_1,0)\int_0^T \bigg(
\Big(\frac{\zeta'}{\zeta}\Big)'(1+iv_1) +
\Big(\frac{t}{2\pi}\Big)^{-iv_1}
\zeta(1+iv_1)\zeta(1-iv_1)A(iv_1) \nonumber \\
&&\qquad\qquad\qquad-\sum_p \frac{\log^2 p} {(p^{1+iv_1}-1)^2}
\bigg)dt\; dv_1\nonumber \\
&&\qquad\qquad+\pi\int_{-T}^T f(0,v_2)\int_0^T \bigg(
\Big(\frac{\zeta'}{\zeta}\Big)'(1+iv_2) +
\Big(\frac{t}{2\pi}\Big)^{-iv_2}
\zeta(1+iv_2)\zeta(1-iv_2)A(iv_2) \nonumber \\
&&\qquad\qquad\qquad-\sum_p \frac{\log^2 p} {(p^{1+iv_2}-1)^2}
\bigg)dt\; dv_2\nonumber \\
&&\qquad\qquad +\pi^2 \int_0^T f(0,0) \log \tfrac{t}{2\pi} dt.
\end{eqnarray}
We note that the final term comprises information on the one-point
correlation function, whereas the second and third integrals in
(\ref{eq:sq1b}) will form a contribution to the two-point
correlation function, as will be discussed at the end of this
section.

In the meantime we will perform sample calculations on three other
terms from (\ref{eq:triplePV}), leaving the remainder of the terms
to the reader, since they all follow a similar pattern.

We consider now
\begin{eqnarray}
\label{eq:sq2}
\int_{-T+i(1-a-b)}^{T+i(1-a-b)}\int_{-T+i(1-a-b)}^{T+i(1-a-b)}
f(v_1,v_2)  I_1(0;iv_2) dv_1dv_2.
\end{eqnarray}
Examining expression (\ref{eq:I1ab}) for $I_1(0;iv_2)$, we see
that the $\big(\frac{\zeta'}{\zeta}\big)'$ term has a second order
pole at $v_2=0$ with residue zero, the term containing
$\zeta(1+iv_2)\zeta(1-iv_2)$ has residue
\begin{equation}
\label{eq:resPV} \int_0^T -i\log^2\tfrac{t}{2\pi} dt ,
\end{equation}
and $B(0,v_2)$ is analytic near $v_2=0$.  The integrand is
analytic in $v_1$, so moving that contour onto the real axis does
not incur any polar contribution.  Thus (\ref{eq:sq2}) equals
\begin{eqnarray}
\int_{-T}^T\int_{-T}^T f(v_1,v_2) I_1(0;v_2) dv_1dv_2 +\pi
\int_{-T}^T f(v_1,0)\int_0^T \log^2 \tfrac{t}{2\pi} dt\; dv_1.
\end{eqnarray}

The next term that deserves consideration is
\begin{eqnarray}
\label{eq:sq4} \int_{-T-i(1-a-b)}^{T-i(1-a-b)} \int_{-T}^T
f(v_1,v_2) I(iv_1,0;-iv_2) dv_1dv_2,
\end{eqnarray}
where again we have used (\ref{eq:switch}) to exchange the
arguments of $I$ preceding the semicolon with those following it.
We note that in these integrals where two consecutive principal
value calculations are made the order of the integrals cannot be
exchanged after the first principal value integral is obtained, so
we always address the principal values starting with the outermost
integral and working inwards.

An inspection of (\ref{eq:Iaab}) reveals that $I(iv_1,0;-iv_2)$ is
not singular at $v_1=0$ (the poles of two $\frac{\zeta'}{\zeta}$
terms cancel), but as we move the $v_2$ contour onto the real axis
we encounter extra difficulties when $v_1=0$, so we will start by
temporarily shifting the $v_1$ contour so that it runs just below
the real axis; we choose below rather than above the axis so as to
avoid the pole at $v_1=v_2$.

Now that the $v_1$ contour does not pass through zero, we can move
the $v_2$ contour onto the real axis and pick up exactly the same
polar contribution as when evaluating (\ref{eq:sq1}).  (An extra
minus sign in the residue compensates for the fact that this time
we need it multiplied by $-i\pi$ since we are half-circling the
pole in the clockwise direction due to the original contour
passing above the real axis.)  Thus, for some $\epsilon>0$,
(\ref{eq:sq4}) is
\begin{eqnarray}
\label{eq:sq4a} &&  \int_{-T}^T \int_{-T-i\epsilon}^{T-i\epsilon}
f(v_1,v_2)I(iv_1,0;-iv_2)dv_1dv_2 \nonumber \\
&&\qquad\qquad+\pi\int_{-T-i\epsilon}^{T-i\epsilon}\int_0^T \bigg(
\Big(\frac{\zeta'}{\zeta}\Big)'(1+iv_1) +
\Big(\frac{t}{2\pi}\Big)^{-iv_1}
\zeta(1+iv_1)\zeta(1-iv_1)A(iv_1) \nonumber \\
&&\qquad\qquad\qquad-\sum_p \frac{\log^2 p} {(p^{1+iv_1}-1)^2}
\bigg) f(v_1,0)dt \;dv_1,
\end{eqnarray}
where the $v_2$ integral is understood as a principal value
integral.

As we move the $v_1$ contour back to the real axis, we encounter a
pole of the first integral above at $v_1=v_2$ (the apparent pole
at $v_1=0$ in fact cancels), and a pole of the second integral at
$v_1=0$ with residue $-i\log \tfrac{t}{2\pi}$, as before.  The
final result is that (\ref{eq:sq4}) equals
\begin{eqnarray}
\label{eq:sq4b} &&\int_{-T}^{T} \int_{-T}^T
f(v_1,v_2)I(iv_1,0;-iv_2)dv_1dv_2 \nonumber \\
&&\qquad\qquad+\pi\int_{-T}^T f(v_1,0)\int_0^T \bigg(
\Big(\frac{\zeta'}{\zeta}\Big)'(1+iv_1) +
\Big(\frac{t}{2\pi}\Big)^{-iv_1}
\zeta(1+iv_1)\zeta(1-iv_1)A(iv_1) \nonumber \\
&&\qquad\qquad\qquad-\sum_p \frac{\log^2 p} {(p^{1+iv_1}-1)^2}
\bigg)dt dv_1\nonumber \\
&&\qquad\qquad+\pi\int_{-T}^T f(v_2,v_2)\int_0^T \bigg(
\Big(\frac{\zeta'}{\zeta}\Big)'(1-iv_2) +
\Big(\frac{t}{2\pi}\Big)^{iv_2}
\zeta(1+iv_2)\zeta(1-iv_2)A(-iv_2) \nonumber \\
&&\qquad\qquad\qquad-\sum_p \frac{\log^2 p} {(p^{1-iv_2}-1)^2}
\bigg)dt dv_2\nonumber \\
&&\qquad\qquad +\pi^2 \int_0^T f(0,0) \log \tfrac{t}{2\pi} dt.
\end{eqnarray}

The last term we will consider (all the others follow in an
identical manner to one of those discussed) is
\begin{eqnarray}
\label{eq:sq7} \int_{-T}^T\int_{-T-i(1-a-b)}^{T-i(1-a-b)}
f(v_1,v_2) I(iv_2,0;-iv_1) dv_1dv_2,
\end{eqnarray}
where as usual we have used (\ref{eq:switch}) to exchange the
arguments of $I$ preceding the semicolon with those following it.
We write this as
\begin{eqnarray}
\label{eq:sq7a} \lim_{\epsilon\rightarrow 0^+}\int_{[-T,T] \atop
|v_2|>\epsilon}\int_{-T-i(1-a-b)}^{T-i(1-a-b)} f(v_1,v_2)
I(iv_2,0;-iv_1) dv_1dv_2.
\end{eqnarray}
This has not changed a thing because the outer integral is
perfectly well-behaved at $v_2=0$, but it means that as we move
the inner integral onto the real axis we avoid the complications
that arise if $v_2=0$.

The inner integral has poles at $v_1=0$ and $v_1=v_2$ with the
usual and now-familiar residues, so (\ref{eq:sq7}) equals
\begin{eqnarray}
\label{eq:sq7b} &&\int_{-T}^{T} \int_{-T}^T
f(v_1,v_2)I(iv_2,0;-iv_1)dv_1dv_2 \nonumber \\
&&\qquad\qquad+\pi\int_{-T}^T f(0,v_2)\int_0^T \bigg(
\Big(\frac{\zeta'}{\zeta}\Big)'(1+iv_2) +
\Big(\frac{t}{2\pi}\Big)^{-iv_2}
\zeta(1+iv_2)\zeta(1-iv_2)A(iv_2) \nonumber \\
&&\qquad\qquad\qquad-\sum_p \frac{\log^2 p} {(p^{1+iv_2}-1)^2}
\bigg)dt \;dv_2\nonumber \\
&&\qquad\qquad+\pi\int_{-T}^T f(v_2,v_2)\int_0^T \bigg(
\Big(\frac{\zeta'}{\zeta}\Big)'(1-iv_2) +
\Big(\frac{t}{2\pi}\Big)^{iv_2}
\zeta(1+iv_2)\zeta(1-iv_2)A(-iv_2) \nonumber \\
&&\qquad\qquad\qquad-\sum_p \frac{\log^2 p} {(p^{1-iv_2}-1)^2}
\bigg)dt\; dv_2,
\end{eqnarray}
where the integrals in $v_1$ and $v_2$ are to be interpreted as
principal value integrals.

All of the other terms in (\ref{eq:triplePV}) can be handled in
exactly the same way as the one of the four treated here.  The
complete result is
\begin{eqnarray}
\label{eq:triplePV1} &&\sum_{0\leq \gamma_1, \gamma_2,\gamma_3<T}
f(\gamma_1-\gamma_2,\gamma_1-\gamma_3)=\frac{1}{(2\pi)^3}\int_{-T}^T\int_{-T}^T
f(v_1,v_2)\nonumber \\
&& \qquad \qquad\times\Bigg(\int_0^T \log ^3 \frac{u}{2\pi} du +
I(iv_1,iv_2;0) + I(0,iv_1;-iv_2) + I(0,iv_2;-iv_1)\nonumber
\\
&&\qquad\qquad\qquad\qquad +
I(-iv_1,-iv_2;0)+I(0,-iv_2;iv_1)+I(0,-iv_1;iv_2)\nonumber
\\&&\qquad\qquad\qquad+I_1(0;iv_2)+I_1(0;iv_1)+I_1(-iv_2;iv_1) +I_1(-iv_2;0)
+I_1(-iv_1;iv_2)\nonumber
\\
&&\qquad\qquad\qquad\qquad+I_1(-iv_1;0)\Bigg)dv_1dv_2\nonumber \\
&&\qquad+\frac{1}{(2\pi)^2} \Bigg(  \int_{-T}^T
f(0,v_2)s(v_2)dv_2+\int_{-T}^T f(v_1,0) s(v_1)dv_1 + \int_{-T}^T
f(v_1,v_1) s(-v_1) dv_1 \nonumber
\\
&&\qquad\qquad+\int_{-T}^T f(0,v_2)s(-v_2)dv_2 + \int_{-T}^T
f(v_1,0) s(-v_1)
dv_1 + \int_{-T}^T f(v_1,v_1) s(v_1) dv_1\Bigg)\nonumber \\
&& \qquad +\frac{1}{2\pi} \int_0^T f(0,0) \log \tfrac{t}{2\pi} dt,
\end{eqnarray}
where
\begin{eqnarray}
&&s(x)=\int_0^T \frac{1}{2} \log^2\tfrac{t}{2\pi} +
\Big(\frac{\zeta'}{\zeta} \Big)'(1+ix)
+\Big(\frac{t}{2\pi}\Big)^{-ix}\zeta(1+ix)\zeta(1-ix)
A(ix)\nonumber
\\
&&\qquad\qquad\qquad\qquad -\sum_p \frac{\log^2p} {(p^{1+ix}-1)^2}
dt.
\end{eqnarray}
 Note that some terms above can be combined if we include the
natural assumption that $f(x,y)=f(-x,-y)$.  All the integrals on
the interval $(-T,T)$ should be considered as principal value
integrals near any poles at the origin or at $v_1=v_2$.

Now we want to identify the terms in the last three lines of
(\ref{eq:triplePV1}) as lower order correlations.  To do this,
note that we can rewrite the triple sum over zeta zeros as
\begin{eqnarray}
\label{eq:zerodecomp} &&\sum_{0<\gamma_1,\gamma_2,\gamma_3<T}
f(\gamma_1-\gamma_2,\gamma_1-\gamma_3) =
\sum_{0<\gamma_1\neq\gamma_2\neq\gamma_3<T}
f(\gamma_1-\gamma_2,\gamma_1-\gamma_3)\nonumber
\\
&&\qquad+\sum_{0<\gamma_1\neq\gamma_2<T} f(\gamma_1-\gamma_2,0) +
\sum_{0<\gamma_1\neq\gamma_3<T} f(0,\gamma_1-\gamma_3)
+\sum_{0<\gamma_1\neq \gamma_2<T}
f(\gamma_1-\gamma_2,\gamma_1-\gamma_2) \nonumber \\
&&\qquad+ \sum_{0<\gamma_1<T} f(0,0).
\end{eqnarray}
The standard result on the density of the zeros of the Riemann
zeta function gives,
\begin{equation}
\label{eq:1p} \sum_{0<\gamma_1<T}1= \frac{T}{2\pi} \log
\frac{T}{2\pi} -\frac{T}{2\pi} +O(1),
\end{equation}
for large $T$, and an expression for the two-point correlation,
derived from the ratios conjecture, is given in
\cite{kn:consna06}:
\begin{eqnarray}
\label{eq:2p}
 &&\sum_{\gamma\neq \gamma'\le T} f(\gamma-\gamma') =\frac{1}{(2\pi)^2}
\int_{-T}^T \int_0^T f(r) \bigg( \log^2 \frac{t}{2\pi} +\left(\frac{\zeta'}{\zeta}\right)'(1+ir)\\
 && \qquad \qquad
+ \left(\frac{t}{2\pi}\right)^ {-ir}
\zeta(1-ir)\zeta(1+ir)A(ir)-B(ir)
+\left(\frac{\zeta'}{\zeta}\right)'(1-ir)\nonumber\\
 && \qquad \qquad
+ \left(\frac{t}{2\pi}\right)^ {ir}
\zeta(1-ir)\zeta(1+ir)A(-ir)-B(-ir) \bigg)~dt ~dr
+O(T^{1/2+\epsilon}),\nonumber
\end{eqnarray}
where $f(z)$ is holomorphic throughout the strip  $|\Im z| <2$, is
real on the real line and satisfies $f(x)\ll 1/(1+x^2)$  as  $x\to
\infty$.  In addition, the integral is to be regarded as a
principal value near $r=0$,
\begin{eqnarray}
A(\eta)=\prod_p\frac{(1-\frac{1}{p^{1+\eta}}) (1-\frac 2 p
+\frac{1}{p^{1+\eta}})}{(1-\frac 1 p )^2},
\end{eqnarray}
and
\begin{eqnarray}
B(\eta)=\sum_p \frac{\log^2 p}{(p^{1+\eta}-1)^2}.
\end{eqnarray}

We see immediately from (\ref{eq:1p}), (\ref{eq:2p}) and
(\ref{eq:zerodecomp}) that
$\sum_{0<\gamma_1\neq\gamma_2\neq\gamma_3<T}
f(\gamma_1-\gamma_2,\gamma_1-\gamma_3)$ is given by the first five
lines of (\ref{eq:triplePV1}), and this is the result presented in
Theorem \ref{theo:triple}.

\subsection{Retrieving the asymptotic result}

\begin{figure}[htbp]
  \begin{center}
    \includegraphics[scale=0.8]
    {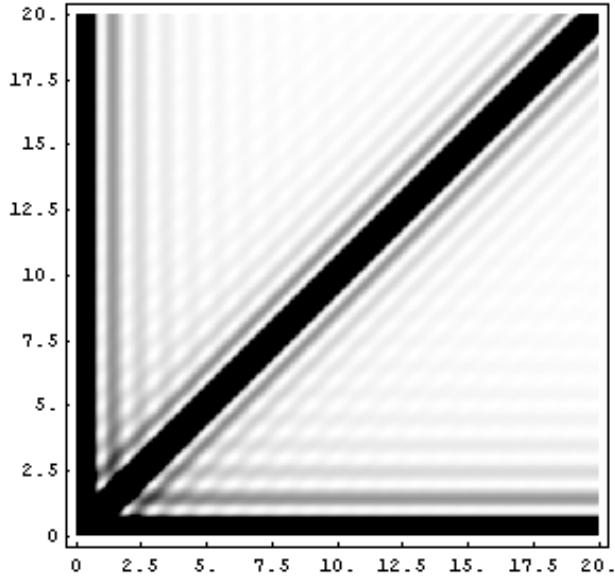}
    \caption{The triple correlation of eigenvalues of random matrices from $U(N)$ in the limit as
    $N\rightarrow\infty$.  That is, we have plotted the $3\times3$ determinant from the
    second line of (\ref{eq:triplimit}).  Note that to compare this picture to Figures \ref{fig:triptheory}
    and \ref{fig:tripnum} the axes of those figures would need to be scaled by the mean density
    of zeros, $(\log \tfrac{T}{2\pi})/(2\pi)$.}
    \label{fig:triprmt}
  \end{center}
\end{figure}

We want to confirm that our formula for the triple correlation of
the Riemann zeros (\ref{eq:triple}) tends  to the limit (plotted
in Figure \ref{fig:triprmt})
\begin{eqnarray}
\label{eq:triplimit} &&\lim_{T\rightarrow
\infty}\frac{1}{\tfrac{T}{2\pi} \log
\tfrac{T}{2\pi}}\sum_{0<\gamma_1\neq \gamma_2\neq \gamma_3<T}
f\bigg(\frac{\log \tfrac{T}{2\pi}}{2\pi}(\gamma_1
-\gamma_2),\frac{\log \tfrac{T}{2\pi}}{2\pi}(\gamma_1-\gamma_3)\bigg) \nonumber\\
&&\qquad= \int_{-\infty}^{\infty} \int
_{-\infty}^{\infty}f(v_1,v_2)\left| \begin{array}{ccc} 1& S(v_1) &
S(v_2)
\\ S(v_1)&1&S(v_1-v_2)\\S(v_2)&S(v_1-v_2)&1\end{array}\right|
dv_1\;dv_2\nonumber\\
&&\qquad=\int_{-\infty}^{\infty} \int
_{-\infty}^{\infty}f(v_1,v_2) \Bigg(
1-\frac{\sin^2\big(\pi(v_1-v_2)\big)} {\pi^2(v_1-v_2)^2}
-\frac{\sin^2\big(\pi v_1\big)} {\pi^2v_1^2} -\frac{\sin^2\big(\pi
v_2\big)} {\pi^2v_2 ^2} \nonumber \\
&&\qquad\qquad \qquad\qquad+2\frac{\sin\big(\pi v_1\big)}{\pi v_1}
\frac{\sin\big(\pi v_2\big)}{\pi v_2} \frac{\sin\big(\pi
(v_1-v_2)\big)}{\pi (v_1-v_2)}\Bigg)dv_1\;dv_2.
\end{eqnarray}

Using (\ref{eq:triple}), with $L=\log \tfrac{T}{2\pi}$, we scale
the variables in the test function by $\tfrac{L}{2\pi}$, and make
a change of variables $v\rightarrow \tfrac{2\pi v}{L}$ in the
$v_1$ and $v_2$ integrals:
\begin{eqnarray}
\label{eq:tripscaled} &&\sum_{0<\gamma_1\neq \gamma_2\neq
\gamma_3<T} f\bigg(\frac{L}{2\pi}(\gamma_1
-\gamma_2),\frac{L}{2\pi}(\gamma_1-\gamma_3)\bigg) \nonumber\\
&&\quad= \frac{1}{L^2 2\pi}
\int_{-\tfrac{TL}{2\pi}}^{\tfrac{TL}{2\pi}}\int_{-\tfrac{TL}{2\pi}}^{\tfrac{TL}{2\pi}}
f(v_1,v_2)\nonumber \\
&& \qquad \qquad\times\Bigg(\int_0^T \log ^3 \frac{u}{2\pi} du +
I\big(\tfrac{2\pi iv_1}{L},\tfrac{2\pi iv_2}{L};0\big) +
I\big(0,\tfrac{2\pi iv_1}{L};-\tfrac{2\pi iv_2}{L}\big) +
I\big(0,\tfrac{2\pi iv_2}{L};-\tfrac{2\pi iv_1}{L}\big)
\\
&&\qquad\qquad\qquad\qquad + I\big(-\tfrac{2\pi
iv_1}{L},-\tfrac{2\pi iv_2}{L};0\big)+I\big(0,-\tfrac{2\pi
iv_2}{L};\tfrac{2\pi iv_1}{L}\big)+I\big(0,-\tfrac{2\pi
iv_1}{L};\tfrac{2\pi iv_2}{L}\big)\nonumber
\\&&\qquad\qquad\qquad+I_1\big(0;\tfrac{2\pi iv_2}{L}\big)+I_1\big(0;\tfrac{2\pi iv_1}{L}\big)+
I_1\big(-\tfrac{2\pi iv_2}{L};\tfrac{2\pi iv_1}{L}\big)
+I_1\big(-\tfrac{2\pi iv_2}{L};0\big) \nonumber
\\
&&\qquad\qquad\qquad\qquad+I_1\big(-\tfrac{2\pi
iv_1}{L};\tfrac{2\pi iv_2}{L}\big)+I_1\big(-\tfrac{2\pi
iv_1}{L};0\big)\Bigg)dv_1dv_2+O(T^{\epsilon}),\nonumber
\end{eqnarray}

Starting with the first term of (\ref{eq:tripscaled}), we see that
asymptotically for large $T$,
\begin{eqnarray}
&&\frac{1}{L^2
2\pi}\int_{-\tfrac{TL}{2\pi}}^{\tfrac{TL}{2\pi}}\int_{-\tfrac{TL}{2\pi}}^{\tfrac{TL}{2\pi}}
f(v_1,v_2)\int_0^T \log^3 \frac{t}{2\pi} dt\; dv_1dv_2\nonumber\\
&& \qquad \sim \frac{T}{2\pi} \log \frac{T}{2\pi}
\int_{-\infty}^{\infty}\int_{-\infty}^{\infty} f(v_1,v_2)
dv_1dv_2.
\end{eqnarray}

A little more work is needed for the other terms of
(\ref{eq:tripscaled}).  Next we consider
\begin{eqnarray}
\label{eq:Iterm1} &&\frac{1}{L^2
2\pi}\int_{-\tfrac{TL}{2\pi}}^{\tfrac{TL}{2\pi}}\int_{-\tfrac{TL}{2\pi}}^{\tfrac{TL}{2\pi}}
f(v_1,v_2)I\big(\tfrac{2\pi i v_1}{L},\tfrac{2\pi i v_2}{L};0\big)
dv_1dv_2.
\end{eqnarray}
Using the definition of $I(\alpha_1,\alpha_2;\beta)$ in
(\ref{eq:Iaab}), we will be a little imprecise and discard any
terms that will not ultimately contribute to the leading-order
$T\log \tfrac{T}{2\pi}$ result.  We thus keep only terms with
third-order poles as $T\rightarrow \infty$, and the leading-order
contribution to (\ref{eq:Iterm1}) is contained in
\begin{eqnarray}
&&\frac{1}{L^2
2\pi}\int_{-\tfrac{TL}{2\pi}}^{\tfrac{TL}{2\pi}}\int_{-\tfrac{TL}{2\pi}}^{\tfrac{TL}{2\pi}}
f(v_1,v_2)\nonumber \\
&&\qquad\qquad \times \int_0^T \bigg(
\Big(\frac{t}{2\pi}\Big)^{-\frac{2\pi i v_1}{L}}
\zeta\big(1-\tfrac{2\pi i v_1}{L}\big) \zeta\big(1+\tfrac{2\pi i
v_1}{L}\big) A\big(\tfrac{2\pi i v_1}{L}\big) \frac{\zeta'}{\zeta}
\big( 1+\tfrac{2\pi i v_2}{L}
-\tfrac{2\pi i v_1}{L}\big) \nonumber \\
&& \qquad\qquad -\Big(\frac{t}{2\pi}\Big)^{-\frac{2\pi i v_1}{L}}
\zeta\big(1-\tfrac{2\pi i v_1}{L}\big) \zeta\big(1+\tfrac{2\pi i
v_1}{L}\big) A\big(\tfrac{2\pi i v_1}{L}\big) \frac{\zeta'}{\zeta}
\big( 1+\tfrac{2\pi i v_2}{L}\big) \nonumber \\
&&\qquad\qquad +\Big(\frac{t}{2\pi}\Big)^{-\frac{2\pi i v_2}{L}}
\zeta\big(1-\tfrac{2\pi i v_2}{L}\big) \zeta\big(1+\tfrac{2\pi i
v_2}{L}\big) A\big(\tfrac{2\pi i v_2}{L}\big) \frac{\zeta'}{\zeta}
\big( 1+\tfrac{2\pi i v_1}{L}
-\tfrac{2\pi i v_2}{L}\big) \nonumber \\
&&\qquad\qquad -\Big(\frac{t}{2\pi}\Big)^{-\frac{2\pi i v_2}{L}}
\zeta\big(1-\tfrac{2\pi i v_2}{L}\big) \zeta\big(1+\tfrac{2\pi i
v_2}{L}\big) A\big(\tfrac{2\pi i v_2}{L}\big) \frac{\zeta'}{\zeta}
\big( 1+\tfrac{2\pi i v_1}{L}\big) \bigg) dt\; dv_1dv_2.
\end{eqnarray}
Keeping only the polar terms and performing the integral over $T$
of the form
\begin{equation}
\int_0^T \Big(\frac{t}{2\pi}\Big)^{-\frac{2 \pi i v}{L}} dt\sim T
e^{-2\pi i v},
\end{equation}
we find (\ref{eq:Iterm1}) is asymptotic to
\begin{eqnarray}
\label{eq:Iterm1a} &&\frac{T}{2\pi} L \int_{-\infty}^{\infty}
\int_{-\infty}^{\infty} f(v_1,v_2) \bigg( -\frac{e^{-2\pi i v_1}}
{i(2\pi)^3 v_1^2(v_2-v_1)} + \frac{e^{-2\pi i v_1}} {i(2\pi)^3
v_1^2 v_2} \nonumber \\
&&\qquad\qquad\qquad- \frac{e^{-2\pi i v_2}} { i (2\pi)^3 v_2^2
(v_1-v_2)} + \frac{e^{-2\pi i v_2}} { i (2\pi)^3 v_1 v_2^2} \bigg)
dv_1 dv_2.
\end{eqnarray}

A similar calculation can be done for
\begin{eqnarray}
\label{eq:Iterm2} &&\frac{1}{L^2
2\pi}\int_{-\tfrac{TL}{2\pi}}^{\tfrac{TL}{2\pi}}\int_{-\tfrac{TL}{2\pi}}^{\tfrac{TL}{2\pi}}
f(v_1,v_2)I\big(0,\tfrac{2\pi i v_1}{L};-\tfrac{2\pi i
v_2}{L}\big) dv_1dv_2,
\end{eqnarray}
and it can be seen to be asymptotic to
\begin{eqnarray}
\label{eq:Iterm2a}&& \frac{T}{2\pi} L \int_{-\infty}^{\infty}
\int_{-\infty}^{\infty} f(v_1,v_2) \bigg( -\frac{e^{2\pi i v_2}}
{i(2\pi)^3 v_1v_2^2} + \frac{e^{2\pi i v_2}} {i(2\pi)^3 v_2^2
(v_1-v_2)} \nonumber \\
&&\qquad\qquad\qquad+ \frac{e^{-2\pi i (v_1-v_2)}} { i (2\pi)^3
v_1 (v_1-v_2)^2} - \frac{e^{-2\pi i (v_1-v_2)}} { i (2\pi)^3 v_2
(v_1-v_2)^2} \bigg) dv_1 dv_2.
\end{eqnarray}
All of the other terms containing a variation of
$I(\alpha_1,\alpha_2;\beta)$ can be obtained from
(\ref{eq:Iterm1a}) and (\ref{eq:Iterm2a}) by a simple swapping of
$v_1$ and $v_2$ or by changing the sign of these variables.

After combining like terms, the exponential terms in the
integrands of (\ref{eq:Iterm1a}) and (\ref{eq:Iterm2a}) and the
other similar integrals sum to
\begin{eqnarray}
 &&-\frac{e^{-2\pi i v_1}} {4i \pi^3 v_1^2(v_2-v_1)} + \frac{e^{2\pi i v_1}} {4i \pi^3
 v_1^2(v_2-v_1)}-\frac{e^{2\pi i v_1}} {4i \pi^3 v_1^2v_2} +\frac{e^{-2\pi i v_1}} {4i \pi^3
 v_1^2v_2}\nonumber \\
 &&\quad-\frac{e^{-2\pi i v_2}} {4i \pi^3 v_2^2(v_1-v_2)}+ \frac{e^{2\pi i v_2}} {4i \pi^3 v_2^2(v_1-v_2)} +\frac{e^{-2\pi i v_2}} {4i \pi^3
 v_1v_2^2} -\frac{e^{2\pi i v_2}} {4i \pi^3
 v_1v_2^2}\nonumber \\
 &&\quad\quad- \frac{e^{2\pi i (v_1-v_2)}} {4i \pi^3 v_1(v_1-v_2)^2}
 + \frac{e^{-2\pi i (v_1-v_2)}} {4i \pi^3 v_1(v_1-v_2)^2} - \frac{e^{-2\pi i (v_1-v_2)}}
 {4i \pi^3 v_2(v_1-v_2)^2} + \frac{e^{2\pi i (v_1-v_2)}} {4i \pi^3
 v_2(v_1-v_2)^2}\nonumber \\
 &&=\frac{\sin(2\pi v_1)} {2\pi^3 v_1^2(v_2-v_1)} - \frac{\sin(2\pi v_1)} {2\pi^3
 v_1^2v_2}+ \frac{\sin(2\pi v_2)} {2\pi^3 v_2^2(v_1-v_2)}
 \nonumber \\
 &&\qquad - \frac{\sin(2\pi v_2)} {2\pi^3 v_1v_2^2} - \frac{\sin\big(2\pi (v_1-v_2)\big)}
 {2\pi^3 v_1(v_1-v_2)^2} + \frac{\sin\big(2\pi (v_1-v_2)\big)} {2\pi^3
 v_2(v_1-v_2)^2}.
 \end{eqnarray}
 Taking a common denominator, we arrive at
 \begin{equation}
 \frac{-\sin(2\pi v_1) +\sin(2\pi v_2) +\sin\big(2\pi
 (v_1-v_2)\big)} {2\pi^3 v_1 v_2(v_1-v_2)},
 \end{equation}
 and use the identity
 \begin{equation}
 \sin x\;\sin y \; \sin(x-y) = \frac{1}{4} \big(\sin(2x-2y)- \sin(2x)
 +\sin(2y)\big)
 \end{equation}
 to show that all the terms in (\ref{eq:tripscaled}) containing
 $I(\alpha_1,\alpha_2;\beta)$ sum to
 \begin{equation}
\frac{T}{2\pi} \log \frac{T}{2\pi}
\int_{-\infty}^{\infty}\int_{-\infty}^{\infty} 2f(v_1,v_2)
\frac{\sin\big(\pi v_1\big)}{\pi v_1} \frac{\sin\big(\pi
v_2\big)}{\pi v_2} \frac{\sin\big(\pi (v_1-v_2)\big)}{\pi
(v_1-v_2)}dv_1\;dv_2.
\end{equation}

The $I_1(\alpha;\beta)$ terms are simpler.  We start with
\begin{eqnarray}\label{eq:I1term}
&&\frac{1}{L^2
2\pi}\int_{-\tfrac{TL}{2\pi}}^{\tfrac{TL}{2\pi}}\int_{-\tfrac{TL}{2\pi}}^{\tfrac{TL}{2\pi}}
f(v_1,v_2)I_1\big(0;\tfrac{2\pi i v_2}{L}\big) dv_1dv_2.
\end{eqnarray}
Picking out just the relevant terms from (\ref{eq:I1ab}), the
leading-order contribution to (\ref{eq:I1term}) is contained in
\begin{eqnarray}
&&\frac{1}{L^2
2\pi}\int_{-\tfrac{TL}{2\pi}}^{\tfrac{TL}{2\pi}}\int_{-\tfrac{TL}{2\pi}}^{\tfrac{TL}{2\pi}}
f(v_1,v_2) \int_0^T \log \frac{t}{2\pi} \bigg(
\Big(\frac{\zeta'}{\zeta}\Big) '\big(1+ \tfrac{2\pi i
v_2}{L}\big)\nonumber \\
&&\qquad\qquad +\Big(\frac{t}{2\pi}\Big)^{-\frac{2\pi i v_2}{L}}
\zeta\big(1+\tfrac{2\pi i v_2}{L}\big) \zeta\big(1-\tfrac{2\pi i
v_2}{L}\big) A\big( \tfrac{2\pi i v_2}{L}\big) \bigg)dt\;dv_1dv_2.
\end{eqnarray}
Expanding the zeta functions around their pole and keeping just
the leading-order term, we then perform the integration over $T$
and arrive at the asymptotic result for (\ref{eq:I1term}) for
large $T$
\begin{eqnarray}
-\frac{T}{2\pi} \log \frac{T}{2\pi} \int_{-\infty}^{\infty}
\int_{-\infty}^{\infty} f(v_1,v_2) \bigg( \frac{1}{(2\pi v_2)^2} -
\frac{e^{-2\pi i v_2}} {(2\pi v_2)^2}\bigg) dv_1dv_2.
\end{eqnarray}

We combine this with
\begin{eqnarray}
&&\frac{1}{L^2
2\pi}\int_{-\tfrac{TL}{2\pi}}^{\tfrac{TL}{2\pi}}\int_{-\tfrac{TL}{2\pi}}^{\tfrac{TL}{2\pi}}
f(v_1,v_2)I_1\big(-\tfrac{2\pi i v_2}{L};0\big) dv_1dv_2\nonumber
\\
&&\quad \sim -\frac{T}{2\pi} \log \frac{T}{2\pi}
\int_{-\infty}^{\infty} \int_{-\infty}^{\infty} f(v_1,v_2) \bigg(
\frac{1}{(2\pi v_2)^2} - \frac{e^{2\pi i v_2}} {(2\pi
v_2)^2}\bigg) dv_1dv_2,
\end{eqnarray}
and the sum gives
\begin{eqnarray}
&&-\frac{T}{2\pi} \log \frac{T}{2\pi} \int_{-\infty}^{\infty}
\int_{-\infty}^{\infty} f(v_1,v_2)  \bigg( \frac{1}{2(\pi v_2)^2}
- \frac{\cos(2\pi v_2)} {2(\pi v_2)^2}\bigg) dv_1dv_2\nonumber \\
&&\quad =-\frac{T}{2\pi} \log \frac{T}{2\pi}
\int_{-\infty}^{\infty} \int_{-\infty}^{\infty} f(v_1,v_2)
\frac{\sin^2(\pi v_2)}{(\pi v_2)^2} dv_1dv_2.
\end{eqnarray}

Since $I_1(\alpha;\beta)$ only depends on the sum of its two
arguments, we see immediately that in the $T\rightarrow \infty$
limit the $I_1$ terms are responsible for
$-\frac{\sin^2\big(\pi(v_1-v_2)\big)} {\pi^2(v_1-v_2)^2}
-\frac{\sin^2\big(\pi v_1\big)} {\pi^2v_1^2} -\frac{\sin^2\big(\pi
v_2\big)} {\pi^2v_2 ^2}$  in (\ref{eq:triplimit}).

\section{Random Matrix Theory}
\label{sect:rmt}

We now use a very similar method to that in Section \ref{sect:rzf}
to derive the triple correlation of eigenvalues of random unitary
matrices.  Of course, there are more elegant methods to do this in
random matrix theory (see Section \ref{sect:gaudin}) but the point
of Section \ref{sect:rmtratios} is that it helps to illuminate the
preceding calculation of the triple correlation of the Riemann
zeros.

If $X$ is an $N\times N$ matrix with complex entries $X=(x_{jk})$,
we let $X^*$ be its conjugate transpose, i.e. $X^*=(y_{jk})$ where
$y_{jk}=\overline{x_{kj}}.$ $X$ is said to be unitary if $XX^*=I$.
We let $U(N)$ denote the group of all $N\times N$ unitary
matrices. This is a compact Lie group and has a Haar measure which
allows us to do analysis.

All of the eigenvalues of $X\in U(N)$ have absolute value 1; we
write them as
  \begin{eqnarray}e^{i\theta_1}, e^{i\theta_2}, \dots , e^{i\theta_N}.\end{eqnarray}

For any sequence of $N$ points on the unit circle there are
matrices in $U(N)$ with these points as eigenvalues. The
collection of all matrices with the same set of eigenvalues
constitutes a conjugacy class in $U(N)$.  Thus, the set of
conjugacy classes can identified with the collection of sequences
of $N$ points on the unit circle.

Weyl's formula asserts that for a function $f:U(N)\to \mathbf C$
which is constant on conjugacy classes,
\begin{eqnarray}\int_{U(N)} f(X) ~d\mbox{Haar}=\int_{[0,2\pi]^N}
f(\theta_1,\dots,\theta_N)dX_N,\end{eqnarray} where
  \begin{eqnarray} dX_N &=&\prod_{1\le j<k\le
N}\big|e^{i\theta_k}-e^{i\theta_j}\big|^2 ~\frac{d\theta_1 \dots
d\theta_N}{N! (2\pi)^N}.\end{eqnarray} Since $N$ will be fixed in
this paper, we will usually write $dX$ in place of $dX_N$.

The characteristic polynomial of a matrix $X$ is denoted
$\Lambda_X(s)$ and is defined by
\begin{eqnarray}\Lambda_X(s)=\det(I-sX^*)=\prod_{n=1}^N(1-se^{-i\theta_n}).
\end{eqnarray}
The roots of $\Lambda_X(s)$ are the eigenvalues of $X$. The
characteristic polynomial  satisfies the functional equation
\begin{eqnarray} \Lambda_X(s)&=&(-s)^N\prod_{n=1}^N e^{-i\theta_n}\prod_{n=1}^N
(1-e^{i\theta_n}/s)\nonumber \\
&=&(-1)^N \det X^* ~s^N~\Lambda_{X^*}(1/s).\end{eqnarray} Note
that
\begin{eqnarray} \label{eqn:fe}
s\frac{\Lambda_X'}{\Lambda_X}(s)+\frac 1 s
\frac{\Lambda_{X^*}'}{\Lambda_{X^*}}\big(\frac 1s\big)=N.
\end{eqnarray}
These characteristic polynomials have value distributions similar
to that of the Riemann zeta-function and form the basis of random
matrix models which predict behavior for the Riemann zeta-function
based on what can be proven about $\Lambda$. Some care has to be
taken in making these comparisons because we are used to thinking
about the zeta-function in a half-plane whereas the characteristic
polynomials are naturally studied on a circle. The translation is
that the 1/2-line corresponds to the unit circle; the half-plane
to the right of the 1/2-line corresponds to the inside of the unit
circle. Note that $\Lambda_X(0)=1$ is the analogue of
$\lim_{\sigma\to \infty}\zeta(\sigma+it)=1$.

We let
\begin{eqnarray}\label{eq:zfunction}
z(x)=\frac{1}{1-e^{-x}}.
\end{eqnarray}
In our formulas for averages of characteristic polynomials the
function $z(x)$ plays the role for random matrix theory that
$\zeta(1+x)$ plays in the theory of moments of the Riemann
zeta-function.

\subsection{Triple correlation by Gaudin's Lemma}
\label{sect:gaudin}
   Let $f(x,y,z)$ be a smooth  function which is periodic with period $2\pi$ in each variable.
We want a formula for
\begin{eqnarray}
T_3(f):=\int_{U(N)}\sideset{}{^*}\sum_{j_1,j_2,j_3}f(\theta_{j_1},\theta_{j_2},\theta_{j_3})dX,
\end{eqnarray}
where the sum is for distinct indices $j_1,j_2,j_3$. It is a
standard result in random matrix theory that by Gaudin's Lemma
(see, for example, \cite{kn:mehta} or \cite{kn:conrey04}), we have
\begin{eqnarray}
T_3(f) \label{eqn:gaudin} = \frac{1}{(2\pi)^3} \int_{[0,2\pi]^3}
f(\theta_1,\theta_2,\theta_3)
  \det_{3\times 3} S_N(\theta_k-\theta_j)~d\theta_3 ~d\theta_2
  ~d\theta_1
\end{eqnarray}
where
\begin{eqnarray}
S_N(\theta)=\frac{\sin \frac {N\theta}{2}}{\sin \frac \theta 2 }.
\end{eqnarray}

\subsection{Triple correlation via the ratios theorem}
\label{sect:rmtratios} We now produce an alternate method of
calculation of the triple correlation for eigenvalues of unitary
matrices. This method mirrors that produced earlier in the paper
for the Riemann zeros, but many steps are cleaner and more obvious
in the random matrix case, not to mention the fact that they are
all rigorous. Therefore the calculation in this section serves to
clarify and support the previous calculation of the triple
correlation of the Riemann zeros.

Let
\begin{eqnarray}
g(z):=\Lambda_X(e^z)=\prod_{j=1}^N\left(1-e^ze^{-i\theta_j}\right).
\end{eqnarray}
Then, since $g(z)$ has zeros at $z=i\theta_j+2\pi i m$,
$m\in\mathbb Z$, by Cauchy's theorem we have for an arbitrary
holomorphic, $2\pi i$ periodic function $f$,
\begin{eqnarray}
\sum_{j=1}^Nf(\theta_j)=\frac{1}{2\pi i}\int_{\mathcal
C}\frac{g'}{g}(z)f(z/i)~dz =\frac{1}{2\pi i}\int_{\mathcal C}e^z
\frac{\Lambda_X'}{\Lambda_X}(e^z)f(z/i)~dz,
\end{eqnarray}
where $\mathcal C$ is a positively oriented contour which encloses
a subinterval of the imaginary axis of length $2\pi$. We choose a
specific path $\mathcal C$ to be the positively oriented rectangle
that has vertices $\delta-\pi i,\delta+\pi i, -\delta+\pi i,
-\delta-\pi i$ where $\delta $ is a small positive number.  Note
that, by periodicity, the integrals on the horizontal segments
cancel each other. Applying this three times, and replacing each
variable by its negative, we have (using the fact that $\mathcal
C$ is unchanged when $z\to-z$),
\begin{eqnarray}M_3(f)&:=&\int_{U(N)}\sum_{j_1=1}^N\sum_{j_2=1}^N
\sum_{j_3=1}^Nf(\theta_{j_1},\theta_{j_2},\theta_{j_3})~dX\\
& =&\frac{-1}{(2\pi i)^3} \int_{\mathcal C}\int_{\mathcal C}
\int_{\mathcal C} e^{-z_1-z_2-z_3}
\int_{U(N)}\frac{\Lambda_X'}{\Lambda_X}(e^{-z_1})
\frac{\Lambda_X'}{\Lambda_X}(e^{-z_2})\frac{\Lambda_X'}{\Lambda_X}(e^{-z_3})~dX\nonumber \\
&& \qquad \times f(iz_1,iz_2,iz_3)~dz_3 ~dz_2 ~dz_1\nonumber
\end{eqnarray}
for a three variable holomorphic periodic function $f$. Notice
that $M_3(f)$ is like $T_3(f)$ except that it is a sum over all
triples $(j_1,j_2,j_3)$ of indices
  between $1$ and $N$ instead of over distinct indices.

Let $\mathcal C_-$ denote the path along the left side of
$\mathcal C$ from $-\delta+\pi i$ down to $-\delta-\pi i$ and  let
$\mathcal C_+$ denote the path along the right side of $\mathcal
C$ from $\delta-\pi i$ up to $\delta+\pi i$.  Ignoring the
integrals over the horizontal paths (because their contribution is
0) we take each variable $z_j$ to be on one or the other of the
two vertical paths $\mathcal C_-$  or $\mathcal C_+$ . In this way
our expression can be written as a sum of eight terms, each term
being a triple integral with each integral on a vertical line
segment either $\mathcal C_-$ or $\mathcal C_+.$ These eight
integrals are analogous to $J_1,\ldots,J_8$ in Section
\ref{sect:contour}.  As we did for the Riemann zeta function, for
each variable $z_j$ which is on $\mathcal C_-$ we use the
functional equation (\ref{eqn:fe}) to replace $
e^{-z_j}\frac{\Lambda_X'}{\Lambda_X}(e^{-z_j}) $  by
$N-e^{z_j}\frac{\Lambda_{X^*}'}{\Lambda_{X^*}}(e^{z_j}). $ In this
way we find (using $X^{-1}=X^*$) that
\begin{eqnarray}&&  \label{eqn:M3}
M_3(f)= \frac{-1}{(2\pi i)^3}\sum_{\epsilon_j\in\{-1,+1\}\atop
j=1,2,3}\int_{\mathcal C_{\epsilon_1}}\int_{\mathcal
C_{\epsilon_2}} \int_{\mathcal C_{\epsilon_3}} \int_{U(N)}
\prod_{j=1}^3 \left(\frac{1-\epsilon_j}{2}N+\epsilon_j
e^{-\epsilon_j z_j}\frac{\Lambda_{X^{\epsilon_j}}'}
{\Lambda_{X^{\epsilon_j}}}
(e^{-\epsilon_j z_j})\right)~dX\\
&&\qquad \qquad \times  f(iz_1,iz_2,iz_3)~dz_3 ~dz_2 ~dz_1.
\nonumber
\end{eqnarray}

We next compute the averages over   $X\in U(N)$ which appear in
the above equation.

We need two instances of the ratios theorem (see \cite{kn:cfz1} or
\cite{kn:cfs05} for statements and proofs of this theorem). Recall
the definition of $z(x)$ from (\ref{eq:zfunction}).

(A) Let $\Re \gamma,\Re \delta>0$.
  \begin{eqnarray}
&& \int_{U(N)}\frac{  \Lambda_X(e^{-\alpha})
\Lambda_{X^*}(e^{-\beta})} { \Lambda_X(e^{-\gamma})
  \Lambda_{X^*}(e^{-\delta})}dX=
   \frac{  z(\alpha+\beta)   z(\gamma+\delta)}
{ z(\alpha+\delta) z(\beta+\gamma) }+e^{-N(\alpha+\beta)}\frac{
z(-\beta-\alpha)   z(\gamma+\delta)} { z(-\beta+\delta)
z(-\alpha+\gamma) }.
\end{eqnarray}

  (B) Let   $\Re \gamma_1,\Re \gamma_2,\Re \delta>0$. Then
  \begin{eqnarray}
&& \int_{U(N)}\frac{
\Lambda_X(e^{-\alpha_1})\Lambda_X(e^{-\alpha_2})
\Lambda_{X^*}(e^{-\beta})} {
\Lambda_X(e^{-\gamma_1})\Lambda_X(e^{-\gamma_2})
  \Lambda_{X^*}(e^{-\delta})}dX=
   \frac{  z(\alpha_1+\beta)  z(\alpha_2+\beta) z(\gamma_1+\delta)z(\gamma_2+\delta)}
{ z(\alpha_1+\delta)z(\alpha_2+\delta) z(\beta+\gamma_1)z(\beta+\gamma_2) }\nonumber\\
&& \qquad +e^{-N(\alpha_1+\beta)}\frac{  z(-\beta-\alpha_1)
z(\alpha_2-\alpha_1)  z(\gamma_1+\delta)z(\gamma_2+\delta)}
{ z(-\beta+\delta)z(\alpha_2+\delta) z(-\alpha_1+\gamma_1)z(-\alpha_1+\gamma_2) }\nonumber\\
&& \qquad +e^{-N(\alpha_2+\beta)}\frac{  z(-\beta-\alpha_2)
z(\alpha_1-\alpha_2)  z(\gamma_1+\delta)z(\gamma_2+\delta)} {
z(-\beta+\delta)z(\alpha_1+\delta)
z(-\alpha_2+\gamma_1)z(-\alpha_2+\gamma_2) }.
\end{eqnarray}

For application to triple correlation, we need   averages of
logarithmic derivatives of the characteristic polynomials.
Differentiating the above formulas leads to the following two
identities, which are the random matrix version of (\ref{eq:I1ab})
and (\ref{eq:Iaab}) respectively, remembering that $z(x)$ plays
the role of $\zeta(1+x)$ and $N$ the role of $\log
\tfrac{t}{2\pi}$.  The products and sums of primes in
(\ref{eq:I1ab}) and (\ref{eq:Iaab}) do not have any counterpart in
random matrix theory.
\begin{proposition} If $\Re \alpha,\Re \beta>0$, then
\begin{eqnarray}
  J(\alpha;\beta):&=&e^{-\alpha-\beta} \int_{U(N)} \frac{\Lambda_X'}{\Lambda_X}(e^{-\alpha})
  \frac{\Lambda_{X^*}'}{\Lambda_{X^*}}(e^{-\beta})  dX\nonumber\\
&=&\left(\frac{z'}{z}\right)'(\alpha+\beta)
+e^{-N(\alpha+\beta)}z(\alpha+\beta)z(-\alpha-\beta).
\end{eqnarray}
\end{proposition}

\begin{proposition}
Let $\Re \alpha_1,\Re \alpha_2,\Re \beta>0$. Then
  \begin{eqnarray}
   J(\alpha_1,\alpha_2;\beta):&=& -e^{-\alpha_1-\alpha_2-\beta}
\int_{U(N)}\frac{\Lambda_X'}{\Lambda_X}(e^{-\alpha_1})\frac{\Lambda_X'}{\Lambda_X}(e^{-\alpha_2})
  \frac{\Lambda_{X^*}'}{\Lambda_{X^*}}(e^{-\beta})  dX \nonumber\\
&=& e^{-N(\alpha_1+\beta)}z(\alpha_1+\beta)z(-\alpha_1-\beta)
\left(\frac{z'}{z}(\alpha_2-\alpha_1)-\frac{z'}{z}(\alpha_2+\beta) \right)\nonumber\\
&&  \qquad +
e^{-N(\alpha_2+\beta)}z(\alpha_2+\beta)z(-\alpha_2-\beta)
\left(\frac{z'}{z}(\alpha_1-\alpha_2)-\frac{z'}{z}(\alpha_1+\beta)
\right) .
\end{eqnarray}
\end{proposition}
Note that the right-hand-side has a simple pole when
$\alpha_1=-\beta$ and when $\alpha_2=-\beta$ but is analytic when
$\alpha_1=\alpha_2$. We also remark that non-constant integrals
with no $X$, or no $X^*$, are 0; for example, if $\Re \alpha, \Re
\beta >0$, then
\begin{eqnarray}
\int_{U(N)}\frac{\Lambda_X'}{\Lambda_X}(e^{-\alpha})~dX=
\int_{U(N)}\frac{\Lambda_{X^*}'}{\Lambda_{X^*}}(e^{-\alpha})~dX=\int_{U(N)}\frac{\Lambda_X'}{\Lambda_X}(e^{-\alpha})
\frac{\Lambda_X'}{\Lambda_X}(e^{-\beta})~dX=0.
\end{eqnarray}

  Thus, letting $dG$ be a shorthand
for $ f(iz_1,iz_2,iz_3)~dz_3~dz_2~dz_1$,
we have, upon expanding (\ref{eqn:M3}),
\begin{eqnarray}&&\label{eqn:M3a}
-(2\pi i)^3M_3(f)=N^3 \int_{\mathcal C_-}\int_{\mathcal
C_-}\int_{\mathcal C_-}~dG+ \int_{\mathcal C_+}\int_{\mathcal
C_+}\int_{\mathcal C_-} J(z_1,z_2;-z_3)~dG
  \\ \nonumber
&&\qquad +\int_{\mathcal C_-}\int_{\mathcal C_+}\int_{\mathcal
C_+} J(z_2,z_3;-z_1)~dG +\int_{\mathcal C_+}\int_{\mathcal
C_-}\int_{\mathcal C_+} J(z_1,z_3;-z_2)~dG
\\ \nonumber
&&\qquad - \int_{\mathcal C_-}\int_{\mathcal C_-}\int_{\mathcal
C_+} ( J(-z_1,-z_2;z_3)+N(
  J(-z_1;z_3)+ J(-z_2;z_3))~dG\\ \nonumber
&&\qquad -\int_{\mathcal C_-}\int_{\mathcal C_+}\int_{\mathcal
C_-} (J(-z_1,-z_3;z_2)+N(J(-z_1;z_2)+J(-z_3;z_2))~dG\\
\nonumber&&\qquad -\int_{\mathcal C_+}\int_{\mathcal
C_-}\int_{\mathcal C_-} (J(-z_2,-z_3;z_1)
+N(J(-z_2;z_1)+J(-z_3;z_1)) ~dG.
\end{eqnarray}
This is very similar to (\ref{eq:triplePV}), except that we have
retained three variables instead of working with the differences
$z_2-z_1$ and $z_3-z_1$.

Now, by the inclusion-exclusion principle,
\begin{eqnarray} \label{eqn:incexc}
T_3(f)&=&M_3(f)-\int_{U(N)}\sum_{j_1,j_2}(f(\theta_{j_1},\theta_{j_1},\theta_{j_2})+
f(\theta_{j_1},\theta_{j_2},\theta_{j_1})+f(\theta_{j_1},\theta_{j_2},\theta_{j_2}))~dX\\
&&\qquad + 2
\int_{U(N)}\sum_{j}f(\theta_{j},\theta_{j},\theta_{j})~dX.
\nonumber
\end{eqnarray}
The pair-correlation sums are evaluated much as above. For
example,
\begin{eqnarray}\label{eqn:pc} &&
\int_{U(N)}\sum_{j_1,j_3}f(\theta_{j_1},\theta_{j_3},\theta_{j_3})dX
=\frac{1}{(2\pi i)^2} \left(-\int_{\mathcal C_-}\int_{\mathcal
C_+}J(-z_1;z_3)
  f(iz_1,iz_3,iz_3)~dz_3~dz_1
\right.\\ \nonumber  &&\qquad \qquad \qquad -\int_{\mathcal
C_+}\int_{\mathcal C_-}J(z_1;-z_3)
  f(iz_1,iz_3,iz_3)~dz_3~dz_1
\\ \nonumber && \qquad \qquad \qquad \qquad \left.
+N^2\int_{\mathcal C_-}\int_{\mathcal C_-}
f(iz_1,iz_3,iz_3)~dz_2~dz_1\right).
  \end{eqnarray}
The one-correlation sum  is
\begin{eqnarray} \label{eqn:onecorr}
&& \int_{U(N)}\sum_{j\le N}f(\theta_j,\theta_j,\theta_j)~dX
=\frac{-N}{2\pi i} \int_{\mathcal C_-}   f(iz,iz,iz)~dz.
\end{eqnarray}

Now we move all of the paths of integration over to the imaginary
axis. When we do this we encounter some poles and have to use
principal value integrals. For example,
\begin{eqnarray}&&
\int_{\mathcal C_+}\int_{\mathcal C_+}\int_{\mathcal C_-}
J(z_1,z_2;-z_3)~dG\\
&&\qquad
= \int_{-\pi i}^{\pi i}\int_{\mathcal C_+}\int_{\mathcal C_-} J(z_1,z_2;-z_3)~dG \nonumber\\
&& \qquad=\int_{-\pi i}^{\pi i}\lim_{\epsilon\to 0^+}\int_{[-\pi
i,\pi i] \atop |z_2-z_1|>\epsilon}
\int_{\mathcal C_-}  J(z_1,z_2;-z_3)~dG \nonumber\\
&& \qquad =-\int_{-\pi i}^{\pi i}\lim_{\epsilon\to 0^+}\int_{[-\pi
i,\pi i] \atop |z_2-z_1|>\epsilon}\lim_{\delta\to 0^+}
\int_{ {[-\pi i,\pi i]\atop |z_3-z_2|>\delta}\atop |z_3-z_1|>\delta} J(z_1,z_2;-z_3)~dG \nonumber\\
&&\qquad \qquad +\pi i \int_{-\pi i}^{\pi i}\lim_{\epsilon\to
0^+}\int_{[-\pi i,\pi i] \atop |z_2-z_1|>\epsilon}
\operatornamewithlimits{Res}_{z_3=z_1}\big(J(z_1,z_2;-z_3)
f(iz_1,iz_2,iz_3)\big)~dz_2~dz_1\nonumber
\\
&&\qquad \qquad +\pi i \int_{-\pi i}^{\pi i}\lim_{\epsilon\to
0^+}\int_{[-\pi i,\pi i] \atop |z_2-z_1|>\epsilon}
\operatornamewithlimits{Res}_{z_3=z_2}\big(J(z_1,z_2;-z_3)
f(iz_1,iz_2,iz_3)\big)~dz_2~dz_1.\nonumber
\end{eqnarray}
The first equality follows because there are no singularities when
we move $z_1$ onto the imaginary axis.  In the second, there are
again no singularities when we move $z_2$ onto the imaginary axis,
but anticipating what comes next we choose to write this integral
with a limit as $\epsilon \to 0^+$.  Finally, the third equality
takes into account the poles at $z_1$ and $z_2$ when we move $z_3$
onto the imaginary axis.  To determine the signs of the three
terms in the final equality, remember the the contour
$\mathcal{C_-}$ is oriented downwards and just to the left of the
imaginary axis.  Thus when it is moved onto the vertical axis it
wraps the singularities in the positive (anti-clockwise)
direction.  Also the residues are multiplied by $\pi i$ instead of
$2\pi i$ because we moved the path onto a path that goes right
through the singularities, i.e. we have not crossed the poles but
rather moved on top of them, and so we have only one-half the
usual contribution of a residual term and the integral in $z_3$
that remains in the first term of the final equality is a
principal value integral.

Now it is easily calculated that
\begin{eqnarray*}
 \operatornamewithlimits{Res}_{z_3=z_1}J(z_1,z_2;-z_3)
=-J(z_2;-z_1)\qquad \mbox{and}\qquad
\operatornamewithlimits{Res}_{z_3=z_2}J(z_1,z_2;-z_3)
=-J(z_1;-z_2).
\end{eqnarray*}
Thus, we have
\begin{eqnarray}   \label{eqn:var1}&&
\int_{\mathcal C_+}\int_{\mathcal C_+}\int_{\mathcal C_-}
J(z_1,z_2;-z_3)~dG \\\nonumber && \qquad = -\int_{-\pi i}^{\pi
i}\lim_{\epsilon\to 0^+}\int_{[-\pi i,\pi i] \atop
|z_2-z_1|>\epsilon}\lim_{\delta\to 0^+}
\int_{ {[-\pi i,\pi i]\atop |z_3-z_2|>\delta}\atop |z_3-z_1|>\delta} J(z_1,z_2;-z_3)~dG \\
&&\qquad \qquad -\pi i \int_{-\pi i}^{\pi i}\lim_{\epsilon\to
0^+}\int_{[-\pi i,\pi i] \atop |z_2-z_1|>\epsilon}  J(z_2;-z_1)
f(iz_1,iz_2,iz_1)~dz_2~dz_1 \nonumber
\\
&&\qquad \qquad -\pi i \int_{-\pi i}^{\pi i}\lim_{\epsilon\to
0^+}\int_{[-\pi i,\pi i] \atop |z_2-z_1|>\epsilon}  J(z_1;-z_2)
f(iz_1,iz_2,iz_2)~dz_2~dz_1.\nonumber
\end{eqnarray}

A slightly different calculation gives
\begin{eqnarray}&&
\int_{\mathcal C_-}\int_{\mathcal C_+}\int_{\mathcal C_+}
J(z_3,z_2;-z_1)~dG\\
&&\qquad
= -\int_{-\pi i}^{\pi i}\int_{\mathcal C_+}\int_{\mathcal C_+} J(z_3,z_2;-z_1)~dG \nonumber\\
&& \qquad=-\int_{-\pi i}^{\pi i}\lim_{\epsilon\to 0^+}\int_{[-\pi
i,\pi i] \atop |z_2-z_1|>\epsilon}
\int_{\mathcal C_+}  J(z_3,z_2;-z_1)~dG\nonumber \\
&&\qquad \qquad + \pi i \int_{-\pi i}^{\pi i}
\operatornamewithlimits{Res}_{z_2=z_1} \bigg(\int_{\mathcal C_+}
J(z_3,z_2;-z_1)~f(iz_1,iz_2,iz_3)~dz_3\bigg) ~dz_1\nonumber\\
&& \qquad=-\int_{-\pi i}^{\pi i}\lim_{\epsilon\to 0^+}\int_{[-\pi
i,\pi i] \atop |z_2-z_1|>\epsilon}
\int_{\mathcal C_+}  J(z_3,z_2;-z_1)~dG \nonumber\\
&&\qquad \qquad + \pi i \int_{-\pi i}^{\pi i}
   \int_{\mathcal C_+}
J(z_3;-z_1)~f(iz_1,iz_1,iz_3)~dz_3 ~dz_1\nonumber
\end{eqnarray}
after moving $z_1$ and $z_2$ onto the imaginary axis.  Now we move
$z_3$ onto the imaginary axis and obtain a residual term from each
of the two integrals above. Thus,
\begin{eqnarray}&&
  \int_{\mathcal C_-}\int_{\mathcal C_+}\int_{\mathcal C_+}
J(z_3,z_2;-z_1)~dG\\
&&\qquad =-\int_{-\pi i}^{\pi i}\lim_{\epsilon\to 0^+}\int_{[-\pi
i,\pi i] \atop |z_2-z_1|>\epsilon}\lim_{\delta\to 0^+} \int_{
{[-\pi i,\pi i]\atop |z_3-z_2|>\delta}\atop |z_3-z_1|>\delta}
J(z_3,z_2;-z_1)~dG\nonumber\\
&& \qquad \qquad -\pi i \int_{-\pi i}^{\pi i}\lim_{\epsilon\to
0^+}\int_{[-\pi i,\pi i] \atop |z_2-z_1|>\epsilon}
\operatornamewithlimits{Res}_{z_3=z_1}
\big(J(z_3,z_2;-z_1)f(iz_1,iz_2,iz_3)\big)~dz_2 ~dz_1\nonumber\\
&& \qquad \qquad -\pi i \int_{-\pi i}^{\pi i}
   \lim_{\epsilon\to 0^+}\int_{[-\pi i,\pi i]\atop |z_3-z_1|>\epsilon}
J(z_3;-z_1)~f(iz_1,iz_1,iz_3)~dz_3 ~dz_1\nonumber\\
&&\qquad \qquad -(\pi i)^2 \int_{-\pi i}^{\pi
i}\operatornamewithlimits{Res}_{z_3=z_1}
\big(J(z_3;-z_1)f(iz_1,iz_1,iz_3)\big)~dz_1.\nonumber
\end{eqnarray}
Now
\begin{eqnarray*}
\operatornamewithlimits{Res}_{z_3=z_1}J(z_3,z_2;-z_1)
=-J(z_2;-z_1) \qquad \mbox{and}\qquad
\operatornamewithlimits{Res}_{z_3=z_1}J(z_3;-z_1) =-N.
\end{eqnarray*}
Thus, we end up with
\begin{eqnarray}\label{eqn:var2}&&
  \int_{\mathcal C_-}\int_{\mathcal C_+}\int_{\mathcal C_+}
J(z_3,z_2;-z_1)~dG\\\nonumber &&\qquad =-\int_{-\pi i}^{\pi
i}\lim_{\epsilon\to 0^+}\int_{[-\pi i,\pi i] \atop
|z_2-z_1|>\epsilon}\lim_{\delta\to 0^+} \int_{ {[-\pi i,\pi
i]\atop |z_3-z_2|>\delta}\atop |z_3-z_1|>\delta}
J(z_3,z_2;-z_1)~dG\\\nonumber && \qquad \qquad -\pi i \int_{-\pi
i}^{\pi i}\lim_{\epsilon\to 0^+}\int_{[-\pi i,\pi i] \atop
|z_2-z_1|>\epsilon} J(z_2;-z_1)f(iz_1,iz_2,iz_1)~dz_2
~dz_1\\\nonumber && \qquad \qquad -\pi i \int_{-\pi i}^{\pi i}
   \lim_{\epsilon\to 0^+}\int_{[-\pi i,\pi i]\atop |z_3-z_1|>\epsilon}
J(z_3;-z_1)~f(iz_1,iz_1,iz_3)~dz_3 ~dz_1\\\nonumber &&\qquad
\qquad -(\pi i)^2 N\int_{-\pi i}^{\pi i}
  f(iz,iz,iz)~dz.
\end{eqnarray}
By a change of variables we have
\begin{eqnarray}\label{eqn:symm}&&
  \int_{\mathcal C_+}\int_{\mathcal C_-}\int_{\mathcal C_-}
J(-z_3,-z_2;z_1)~dG=
  -\int_{\mathcal C_-}\int_{\mathcal C_+}\int_{\mathcal C_+}
J(z_3,z_2;-z_1)~dG.
\end{eqnarray}
We can also make use of the symmetry
$J(\alpha_1,\alpha_2;\beta)=J(\alpha_2,\alpha_1;\beta)$. In this
way we can use (\ref{eqn:var1}) and (\ref{eqn:var2}) to replace
all of the integrals involving a three-variable $J$ in terms of
principal value integrals.

In a similar way we rewrite, for example, the two variable $J$-
integral as
\begin{eqnarray}\label{eqn:twovarJ}&&
\int_{\mathcal C_+}\int_{\mathcal C_-}J(z_1;-z_3)
f(iz_1,iz_3,iz_3)~dz_3~dz_1\\ && \qquad =-  \int_{-\pi i}^{\pi
i}\lim_{\epsilon\to 0^+} \int_{[-\pi i,\pi i]\atop
|z_3-z_1|>\epsilon} J(z_1;-z_3) f(iz_1,iz_3,iz_3)~dz_3~dz_1\nonumber\\
\nonumber&& \qquad \qquad-\pi i N \int_{-\pi i}^{\pi i}
f(iz,iz,iz)~dz.
\end{eqnarray}
Using (\ref{eqn:M3}) -- (\ref{eqn:twovarJ}), after the
substitutions $iz_j=\theta_j$, $j=1,2,3$, and  some
simplification, we have
\begin{eqnarray}&&
T_3(f)  =\frac{1}{(2\pi)^3}\int_{-\pi}^\pi \lim_{\epsilon\to
0^+}\int_{[-\pi,\pi]\atop |\theta_2-\theta_1|>\epsilon
}\lim_{\delta\to 0^+} \int_{{[-\pi,\pi]\atop
|\theta_3-\theta_1|>\delta}\atop |\theta_3-\theta_2|>\delta}
  \bigg(  J(i\theta_1,i\theta_2;-i\theta_3)
  +  J(i\theta_1,i\theta_3;-i\theta_2)\nonumber\\ &&\qquad     \qquad + J(i\theta_2,i\theta_3;-i\theta_1)
     +
  J(-i\theta_1,-i\theta_2;i\theta_3)  +
J(-i\theta_1,-i\theta_3;i\theta_2)   +J(-i\theta_2,-i\theta_3;i\theta_1)\nonumber\\
&& \qquad \qquad \qquad +N\big(
  J(-i\theta_1;i\theta_3)+ J(-i\theta_2;i\theta_3)+J(-i\theta_1;i\theta_2)+
  J(-i\theta_3;i\theta_2)\nonumber\\
&&\qquad \qquad \qquad \qquad +
J(-i\theta_2;i\theta_1)+J(-i\theta_3;i\theta_1)\big)+N^3 \bigg)
f(\theta_1,\theta_2,\theta_3) ~d\theta_1 ~d\theta_2 ~d\theta_3,
\end{eqnarray}
all of the pair--correlation and one-correlation terms having
cancelled.  The principal value integrals with the limits in
$\epsilon$ and $\delta$ are no longer needed because it can be
checked, using Mathematica for example, that the integrand here is
\begin{eqnarray}
=\det_{3\times 3}S_N(\theta_k-\theta_j)f(\theta_1,
\theta_2,\theta_3)
\end{eqnarray}
which is entire.  Now our formula agrees with (\ref{eqn:gaudin}).

This concludes our alternate proof of triple correlation for
unitary matrices.


\end{document}